\documentstyle[12pt]{article}
\input amssym.def
\begin {document}
\def \Z{\bf Z}
\def \C{\bf C}

\def \Q{\Bbb Q}
\def \N{\Bbb N}

\def \bZ{{\bf Z}}
\def \TZ{{1\over T}{\bf Z}}

\def \fg{\frak g}

\def \Res{{\rm Res}}
\def \End{{\rm End}\;}

\def \mod{{\rm mod}\;}

\def \<{\langle} 
\def \>{\rangle} 

\def \a{\alpha }

\def \b{\beta }

\def \be{\begin{equation}\label}
\def \ee{\end{equation}}
\def \bl{\begin{lem}\label}
\def \el{\end{lem}}
\def \bt{\begin{thm}\label}
\def \et{\end{thm}}
\def \bp{\begin{prop}\label}
\def \ep{\end{prop}}
\def \br{\begin{rem}\label}
\def \er{\end{rem}}
\def \bc{\begin{coro}\label}
\def \ec{\end{coro}}
\def \bd{\begin{de}\label}
\def \ed{\end{de}}
\def \pf{{\bf Proof. }}

\newtheorem{thm}{Theorem}[section]
\newtheorem{prop}[thm]{Proposition}
\newtheorem{coro}[thm]{Corollary}

\newtheorem{lem}[thm]{Lemma}
\newtheorem{rem}[thm]{Remark}
\newtheorem{de}[thm]{Definition}

\makeatletter
\@addtoreset{equation}{section}

\makeatother
\makeatletter

\baselineskip=16pt
\begin{center}{\Large \bf Generalized vertex algebras generated by
parafermion-like vertex operators}

\vspace{0.5cm}
Yongcun Gao\\
Department of Mathematics, Nankai University, Tianjin 300071, China\\
Haisheng Li\footnote{Partially supported by NSF grant
DMS-9970496}\\ Department of Mathematical Sciences, Rutgers University,
Camden, NJ 08102
\end{center}

{\bf Abstract} 
It is proved that for a vector space $W$, any set of parafermion-like 
vertex operators on $W$ in a certain canonical way generates a generalized 
vertex algebra in the sense of [DL2] with $W$ as a natural module. 
This result generalizes a result of [Li2]. As an application, generalized 
vertex algebras are constructed from Lepowsky-Wilson's $Z$-algebras of 
any nonzero level.

\section{Introduction}

Vertex operator algebras, introduced in mathematics ([B], [FLM]),
are known essentially to be chiral algebras, 
introduced in physics ([BPZ], [MS]). In terms of physical language, 
chiral algebras for bosonic field theories are vertex operator algebras
while chiral algebras are vertex operator superalgebras
for fermionic field theories.
In physics, further generalizations of bosons and fermions 
are parafermions ([ZF1-2], [G]), where the chiral algebras
were called parafermion algebras.
 Independently (and earlier), $Z$-operators and $Z$-algebras 
([LW1-6], [LP1-2]) were introduced in mathematics
to study standard modules for affine Lie algebras.
In [DL1-3] and [LP1-2], the relations between $Z$-operators and 
parafermion operators were clarified.
Furthermore, in [DL2], Lepowsky-Wilson's $Z$-algebras were put
into larger, more natural algebras where the notions of
generalized vertex (operator) algebra and abelian intertwining algebras
were introduced in [DL2]. Generalized algebraic structures 
associated with rational lattices were also studied in [M]
and a notion called vertex operator para-algebra
was independently introduced and studied in [FFR].

Roughly speaking, parafermion algebras
are generalized vertex operator algebras.
For bosonic or fermionic field operators $a(z)$ and $b(z)$, 
the locality amounts to
\begin{eqnarray}
(z_{1}-z_{2})^{k}a(z_{1})b(z_{2})
=(-1)^{|a||b|}(z_{1}-z_{2})^{k}b(z_{2})a(z_{1})
\end{eqnarray}
for some non-negative integer $k$, where $|a|=0$ if $a(z)$ is bosonic
and $1$ if $a(z)$ is fermionic.
Parafermions are always associated to an abelian group $G$
equipped with a $\C^{\times}$-valued alternating $\Z$-bilinear
form
$c(\cdot,\cdot)$ and a $\C/2\Z$-valued $\Z$-bilinear 
form $(\cdot,\cdot)$ on $G$.
For parafermion operators $a(z), b(z)$ with gradings $g,h\in G$, 
the following relation holds:
\begin{eqnarray}
(z_{1}-z_{2})^{k+(g,h)}a(z_{1})b(z_{2})
=(-1)^{k}c(g,h)(z_{2}-z_{1})^{k+(g,h)}b(z_{2})a(z_{1})
\end{eqnarray}
for some nonnegative integer $k$.

In physical literatures, a chiral algebra is often
described by a set of generating field operators and
a certain set of relations (such as operator product expansions). 
Now we know that 
a certain set of field operators on a vector space
indeed gives rise to a vertex operator (super)algebra.
A result proved in [Li2] is
that any set of mutually local vertex operators on a vector space $W$
generates a canonical vertex superalgebra with $W$ as a module.
(This result was extended in [Li3] for twisted modules.)
This is an analogue of the simple fact in linear algebra that
any set of mutually commutative endomorphisms on a vector space $U$
generates a commutative associative algebra with $U$ as a module.
(See [FKRW], [LZ], [MN], [MP] and [X] for other related interesting results.)

As the main result of this paper, we extend the results of [Li2-3]
for the notion of generalized vertex algebra [DL2]. 
This paper is modeled on [Li2], however, 
the two key theorems (Propositions 3.8 and 3.13) require
essentially new proofs.

It seems that the most general and natural notion is the one of abelian
intertwining algebra [DL2]. An abelian
intertwining algebra by definition is associated to an abelian group $G$
and $\C^{\times}$-valued functions $F$ on $G\times G\times G$ 
and $\Omega$ on $G\times G$ satisfying certain conditions.
It would be nice to extend 
our result to the notion of abelian intertwining algebra.
When we were trying, we found that
the extension is almost straightforward
except for extending the two key theorems (Propositions 3.8 and 3.13)
we need certain identities purely about $F$ and $\Omega$.
We can prove that one of the identities follows from
the assumptions on $F$ and $\Omega$, but we are not be able to
prove the others. We hope to discuss this issue in some other place.

This paper consists of four sections including this introduction as
the first section. In Section 2, we recall some basic definitions
and results from [DL2]. The main result of this paper is given
in Section 3. In Section 4, we construct
canonical generalized vertex algebras from
$Z$-algebras of any nonzero level.

\section{Definition and duality for generalized vertex algebras}
This section is preliminary.
In this section, we recall from [DL2] the basic definitions
(of generalized vertex algebra and module) and basic 
duality properties.

First, let us briefly review some formal variable calculus.
(Best references are [FLM] and [FHL].)
Throughout this paper, $z, z_{0},z_{1},z_{2}$ and $x, y$ will be
mutually commuting (independent) formal variables. 
We shall use $\N$ for the set of all nonnegative integers,
$\Z_{+}$ for the set of positive integers and $\C$ for the set 
of complex numbers. All vector spaces are assumed to be over $\C$.

For a vector space $U$, set
\begin{eqnarray}
U\{z\}=\left\{\sum_{n\in \C}u(n)z^{n}\;|\; u(n)\in U
\;\;\;\mbox{ for }n\in \C\right\}.
\end{eqnarray}
Let  $D=d/dz$ be the formal differential operator on $U\{z\}$:
\begin{eqnarray}
D\left(\sum_{n\in \C}u(n)z^{n}\right)=\sum_{n\in \C}nu(n)z^{n-1}. 
\end{eqnarray}
The formal residue operator $\Res_{z}$ from 
$U\{z\}$ to $U$ is defined by
\begin{eqnarray}
\Res_{z}u(z)=u(-1)
\end{eqnarray}
for $u(z)=\sum_{n\in \C}u(n)z^{n}\in U\{z\}$.

The following are useful subspaces of $U\{z\}$:
\begin{eqnarray}
U[[z,z^{-1}]]&=&\left\{ \sum_{n\in \Z}u(n)z^{n}\;|\; u(n)\in U
\;\;\;\mbox{ for }n\in \Z\right\},\\
U((z))&=&\left\{ \sum_{n\in \Z}u(n)z^{n}\in U[[z,z^{-1}]]\;|\; u(n)=0
\;\;\;\mbox{ for }n\;\;\mbox{sufficiently small}\right\},\\
U[[z]]&=&\left\{ \sum_{n\in \Z}u(n)z^{n}\in U[[z,z^{-1}]]\;|\; u(n)=0
\;\;\;\mbox{ for }n<0\right\}.
\end{eqnarray}
A typical element of $\C[[z,z^{-1}]]$ is the formal Fourier expansion
of the delta-function at $0$:
\begin{eqnarray}
\delta(z)=\sum_{n\in \Z}z^{n}.
\end{eqnarray}
Its fundamental property is:
\begin{eqnarray}
f(z)\delta(z)=f(1)\delta(z)\;\;\;\mbox{ for }\; f(z)\in \C[z,z^{-1}].
\end{eqnarray}
For $\a\in \C$, by definition,
\begin{eqnarray}
(z_{1}-z_{2})^{\alpha}=\sum_{i\ge 0}{\a\choose i}(-1)^{i}z_{1}^{\a-i}z_{2}^{i}.
\end{eqnarray}
Then
\begin{eqnarray}
\delta\left(\frac{z_{1}-z_{2}}{z_{0}}\right)
=\sum_{n\in \Z}\left(\frac{z_{1}-z_{2}}{z_{0}}\right)^{n}
=\sum_{n\in \Z}\sum_{i\ge 0}{n\choose i}(-1)^{i}z_{0}^{-n}z_{1}^{n-i}z_{2}^{i}.
\end{eqnarray}
We have the following fundamental properties of delta function
([FLM], [FHL], [Le], [Zhu]):

\bl{ldeltajacobi}
For $\a\in \C$,
\begin{eqnarray}\label{edeltaalpha}
z_{0}^{-1}\left(\frac{z_{1}-z_{2}}{z_{0}}\right)^{\a}
\delta\left(\frac{z_{1}-z_{2}}{z_{0}}\right)
=z_{1}^{-1}\left(\frac{z_{0}+z_{2}}{z_{1}}\right)^{-\a}
\delta\left(\frac{z_{0}+z_{2}}{z_{2}}\right);
\end{eqnarray}
For $r,s,k\in \Z$ and for $p(z_{1},z_{2})\in \C[[z_{1},z_{2}]]$,
\begin{eqnarray}
& &z_{0}^{-1}\delta\left(\frac{z_{1}-z_{2}}{z_{0}}\right)
z_{1}^{r}z_{2}^{s}(z_{1}-z_{2})^{k}p(z_{1},z_{2})
-z_{0}^{-1}\delta\left(\frac{z_{2}-z_{0}}{-z_{0}}\right)
z_{1}^{r}z_{2}^{s}(-z_{2}+z_{1})^{k}p(z_{1},z_{2})\nonumber\\
&=&z_{2}^{-1}\delta\left(\frac{z_{1}-z_{0}}{z_{2}}\right)
(z_{2}+z_{0})^{r}z_{2}^{s}z_{0}^{k}p(z_{2}+z_{0},z_{2}).
\end{eqnarray}
In particular,
\begin{eqnarray}
z_{0}^{-1}\delta\left(\frac{z_{1}-z_{2}}{z_{0}}\right)-
z_{0}^{-1}\delta\left(\frac{z_{2}-z_{0}}{-z_{0}}\right)
=z_{2}^{-1}\delta\left(\frac{z_{1}-z_{0}}{z_{2}}\right).
\end{eqnarray}
\el

Note that (\ref{edeltaalpha}) is equivalent to
\begin{eqnarray}
(z_{1}-z_{2})^{\a}
z_{0}^{-1}\delta\left(\frac{z_{1}-z_{2}}{z_{0}}\right)
=z_{1}^{-1}z_{0}^{\a}\left(\frac{z_{0}+z_{2}}{z_{1}}\right)^{-\a}
\delta\left(\frac{z_{0}+z_{2}}{z_{2}}\right)\;\;\;\mbox{ for }\a\in \C.
\end{eqnarray}

A generalized vertex algebra by definition is associated to
an abelian group $G$, a symmetric 
 $\C /2{\bf Z}$-valued ${\bf Z}$-bilinear form 
(not necessarily nondegenerate) on $G$:
\begin{eqnarray}
(g,h)\in \C/2{\bf Z}\;\;\;\mbox{ for }g,h\in G
\end{eqnarray}
and $c(\cdot,\cdot)$ is a ${\bf C}^{\times}$-valued 
alternating ${\bf Z}$-bilinear form on $G$.

{\em A generalized vertex algebra} associated
with the group $G$ and the forms $(\cdot,\cdot)$ and $c(\cdot,\cdot)$
is a $G$-graded vector space
\begin{eqnarray}
V=\coprod_{g\in G}V^{g},
\end{eqnarray}
equipped with a linear map
\begin{eqnarray}
Y: & &V \rightarrow (\End V)\{z\}\nonumber\\ 
& &u\mapsto Y(u,z)=\sum_{n\in \C}u_{n}z^{-n-1}
\end{eqnarray}
and with a distinguished vector ${\bf 1}\in V^{0}$,
called the {\em vacuum vector},
satisfying the following conditions for
$g,h\in G,\; u,v\in V$ and $l\in \C$:
\begin{eqnarray}
& &u_{l}V^{h}\subset V^{g+h}\;\;\;\mbox{ if }u\in V^{g};\\
& &u_{l}v=0\;\;\;\mbox{ if the real part of $l$ is sufficiently large};
\label{etruncation}\\
& &Y({\bf 1},z)=1;\label{evacuum}\\
& &Y(v,z){\bf 1}\in V[[z]]\;\;\mbox{ and }
v_{-1}{\bf 1}\;(=\lim_{z\rightarrow 0}Y(v,z){\bf 1})=v;\\
& &Y(u,z)|_{V^{h}}=\sum_{n\equiv (g,h)\mod \bZ}v_{n}z^{-n-1}
\;\;\;\mbox{ if }u\in V^{g}
\end{eqnarray}
(i.e., $n+2\bZ\equiv (g,h)\;\mod \bZ/2\bZ$);
\begin{eqnarray}\label{egjacobi0}
& &z_{0}^{-1}\left(\frac{z_{1}-z_{2}}{z_{0}}\right)^{(g,h)}
\delta\left(\frac{z_{1}-z_{2}}{z_{0}}\right)Y(u,z_{1})Y(v,z_{2})\nonumber\\
& &-c(g,h)z_{0}^{-1}\left(\frac{z_{2}-z_{1}}{z_{0}}\right)^{(g,h)}
\delta\left(\frac{z_{2}-z_{1}}{-z_{0}}\right)Y(v,z_{2})Y(u,z_{1})\nonumber\\
&=&z_{2}^{-1}\delta\left(\frac{z_{1}-z_{0}}{z_{2}}\right)Y(Y(u,z_{0})v,z_{2})
\left(\frac{z_{1}-z_{0}}{z_{2}}\right)^{-g}
\end{eqnarray}
(the {\em generalized Jacobi identity}) if $u\in V^{g},\; v\in V^{h}$, where
\begin{eqnarray}
\delta\left(\frac{z_{1}-z_{0}}{z_{2}}\right)
\left(\frac{z_{1}-z_{0}}{z_{2}}\right)^{-g}\cdot w=
\left(\frac{z_{1}-z_{0}}{z_{2}}\right)^{-(g,g')}
\delta\left(\frac{z_{1}-z_{0}}{z_{2}}\right)w
\end{eqnarray}
for $w\in V^{g'},\; g'\in G$.

This completes the definition. The map $Y$ is called the 
{\em vertex operator map}.
The generalized vertex algebra is denoted by
$$(V,Y, {\bf 1}, G, c(\cdot,\cdot),(\cdot,\cdot))$$ 
or briefly, by  $V$.

\br{rchange}
{\em Note that we here slightly generalize
the original definition in [DL2] where $(\cdot,\cdot)$ was assumed
to be $({1\over T}\Z)/2\Z$-valued, where $T$ is a positive integer
called the level.
The main reason for this generalization is to include 
the generalized vertex algebras
associated to affine algebras with a non-rational level.
On the other hand, if $G$ is finite, 
there exists a positive integer $T$ such that
$(\cdot,\cdot)$ ranges in $\TZ/2\Z$.}
\er

We recall the following remarks from [DL2]:

\br{rvoa}
{\em If $G=0$, the notion of generalized vertex algebra reduces 
to the notion of vertex algebra.
If $G=\Z/2\Z$ with $(m+\Z,n+\Z)=mn+\Z$, 
the notion of generalized vertex algebra reduces 
to the notion of vertex superalgebra, noting that $c(\cdot,\cdot)=1$.}
\er

\br{rgvoa}
{\em A {\em generalized vertex operator algebra} is 
a generalized vertex algebra $V$
associated to a finite group $G$ with $c(\cdot,\cdot)=1$ and
$(\cdot,\cdot)$ being nondegenerate and
furthermore, it is equipped with another distinguished vector 
$\omega\in V_{2}^{0}$, called the Virasoro vector, such that
\begin{eqnarray}
& &[L(m),L(n)]=(m-n)L(m+n)+{1\over 12}(m^{3}-m)\delta_{m+n,0}c
\;\;\;\mbox{ for }m,n\in \Z,\\
& &[L(-1),Y(v,z)]={d\over dz}Y(v,z)\;\;\;\mbox{ for }v\in V,
\end{eqnarray}
where $Y(\omega,z)=\sum_{m\in \Z}L(m)z^{-m-2}$ and $c$ is a 
complex number, called the {\em rank} of $V$, and such that
\begin{eqnarray}
& &V=\coprod_{n\in \C}V_{n}\;\;\mbox{ where }
V_{n}=\{v\in V\:|\;L(0)v=n v\}\;\;\mbox{ for } n\in \C;\\
& &V^{g}=\coprod_{n\in \C}V_{n}^{g}\;\;\;\;
\mbox{(where $V_{n}^{g}=V_{n}\cap V^{g}$)}\;\;\mbox{ for }g\in G;\\
& &\dim V_{n}<\infty \;\;\;\mbox{ for }n\in \C,\\
& &V_{n}=0\;\;\;\mbox{ for $n$ whose real part is sufficiently small}.
\end{eqnarray}}
\er

\bp{pdlduality}
\mbox{[DL2]} In the presence of all the axioms
except the generalized Jacobi identity in defining the notion 
of generalized vertex algebra,
the generalized Jacobi identity is equivalent to the following
generalized weak commutativity and associativity:

(A) For $g_{1},g_{2}\in G$ and $v_{1}\in V^{g_{1}},\; v_{2}\in V^{g_{2}}$,
there exists a nonnegative integer $k$ such that
\begin{eqnarray}
& &(z_{1}-z_{2})^{k+(g_{1},g_{2})}Y(v_{1},z_{1})Y(v_{2},z_{2})\nonumber\\
&=&(-1)^{k}c(g_{1},g_{2})(z_{2}-z_{1})^{k+(g_{1},g_{2})}
Y(v_{2},z_{2})Y(v_{1},z_{1}).
\end{eqnarray}
(B) For $g_{1},g_{2},h\in G$ and 
$v_{1}\in V^{g_{1}},\; v_{2}\in V^{g_{2}},\; w\in V^{h}$,
there exists a nonnegative integer $l$ such that
\begin{eqnarray}
(z_{0}+z_{2})^{l+(g_{1},h)}Y(v_{1},z_{0}+z_{2})Y(v_{2},z_{2})w
=(z_{2}+z_{0})^{l+(g_{1},h)}Y(Y(v_{1},z_{0})v_{2},z_{2})w,
\end{eqnarray}
where $l$ is independent of $v_{2}$.
\ep


The following result was also due to [DL2]:

\bp{pdlgvoa}
In the presence of all the axioms
except the generalized Jacobi identity in defining the notion 
of generalized vertex algebra, the generalized 
Jacobi identity follows from the generalized weak commutativity
and the following property
\begin{eqnarray}\label{ebracket-derivative}
[D,Y(v,z)]={d\over dz}Y(v,z)\;\;\;\mbox{ for }v\in V,
\end{eqnarray}
where $D$ is an endomorphism of $V$ defined by 
\begin{eqnarray}
D(v)=v_{-2}{\bf 1}\;\;\;\mbox{ for }v\in V.
\end{eqnarray}
\ep

\pf We here give a slightly different proof by generalizing the proof
given in [Li2] for the corresponding result for vertex superalgebras.

First, just as in [Li2], from the vacuum property and 
(\ref{ebracket-derivative}) 
we have
\begin{eqnarray}
Y(v,z){\bf 1}=e^{zD}v\;\;\;\mbox{ for }v\in V.
\end{eqnarray}

Second, (\ref{ebracket-derivative}) is equivalent to 
the following conjugation formula:
\begin{eqnarray}\label{econjugation}
e^{z_{0}D}Y(v,z)e^{-z_{0}D}=Y(v,z+z_{0}).
\end{eqnarray}

Third, we shall derive a skew-symmetry. 
Let $u\in V^{g},\; v\in V^{h},\; g,h\in G$.
Then there exists a nonnegative integer $k$ such that
\begin{eqnarray}
(z_{1}-z_{2})^{k+(g,h)}Y(u,z_{1})Y(v,z_{2})=(-1)^{k}
c(g,h)(z_{2}-z_{1})^{k+(g,h)}Y(v,z_{2})Y(u,z_{1})
\end{eqnarray}
and such that
\begin{eqnarray}\label{epositive}
z^{k+(g,h)}Y(v,z)u\in V[[z]].
\end{eqnarray}
Then
\begin{eqnarray}
& &(z_{1}-z_{2})^{k+(g,h)}Y(u,z_{1})Y(v,z_{2}){\bf 1}\nonumber\\
&=&(-1)^{k}c(g,h)(z_{2}-z_{1})^{k+(g,h)}Y(v,z_{2})Y(u,z_{1}){\bf 1}\nonumber\\
&=&(-1)^{k}c(g,h)(z_{2}-z_{1})^{k+(g,h)}Y(v,z_{2})e^{z_{1}D}u\nonumber\\
&=&(-1)^{k}c(g,h)e^{z_{1}D}(z_{2}-z_{1})^{k+(g,h)}Y(v,z_{2}-z_{1})u\nonumber\\
&=&(-1)^{k}c(g,h)e^{z_{1}D}(e^{\pi i}z_{1}+z_{2})^{k+(g,h)}
Y(v,e^{\pi i}z_{1}+z_{2})u.
\end{eqnarray}
We are using (\ref{epositive}). Now, it is safe 
for us to replace $z_{2}$ with $0$.
In this way we get
\begin{eqnarray}
z_{1}^{k+(g,h)}Y(u,z_{1})v
=c(g,h)e^{\pi i(g,h)}z_{1}^{k+(g,h)}e^{z_{1}D}Y(v,e^{\pi i}z_{1})u.
\end{eqnarray}
Then we obtain the following skew-symmetry:
\begin{eqnarray}\label{eskewsymmetry}
Y(u,z_{1})v=c(g,h)e^{\pi i(g,h)}e^{z_{1}D}Y(v,e^{\pi i}z_{1})u.
\end{eqnarray}

Next, we prove the generalized weak associativity.
Let $u\in V^{g_{1}},\; v\in V^{g_{2}},\; w\in V^{g_{3}}$.
Let $k$ be a nonnegative integer (only depending on $u,w$) such that
\begin{eqnarray}
(z_{1}-z_{2})^{k+(g_{1},g_{3})}Y(u,z_{1})Y(w,z_{2})
=(-1)^{k}c(g_{1},g_{3})(z_{2}-z_{1})^{k+(g_{1},g_{3})}Y(w,z_{2})Y(u,z_{1}).
\end{eqnarray}
Then using the skew-symmetry (\ref{eskewsymmetry}) and 
the conjugation formula (\ref{econjugation}) we obtain
the following generalized associativity relation
\begin{eqnarray}
& &(z_{0}+z_{2})^{k+(g_{1},g_{3})}Y(u,z_{0}+z_{2})Y(v,z_{2})w\nonumber\\
&=&c(g_{2},g_{3})e^{\pi i(g_{2},g_{3})}(z_{0}+z_{2})^{k+(g_{1},g_{3})}
Y(u,z_{0}+z_{2})e^{z_{2}D}Y(w,e^{\pi i}z_{2})v\nonumber\\
&=&c(g_{2},g_{3})e^{\pi i(g_{2},g_{3})}(z_{0}+z_{2})^{k+(g_{1},g_{3})}
e^{z_{2}D}Y(u,z_{0})Y(w,e^{\pi i}z_{2})v\nonumber\\
&=&(-1)^{k}c(g_{2},g_{3})c(g_{1},g_{3})e^{\pi i(g_{2},g_{3})}
(e^{\pi i}z_{2}-z_{0})^{k+(g_{1},g_{3})}
e^{z_{2}D}Y(w,e^{\pi i}z_{2})Y(u,z_{0})v\nonumber\\
&=&(-1)^{k}e^{\pi i(g_{2},g_{3})}
(e^{\pi i}z_{2}-z_{0})^{k+(g_{1},g_{3})}
c(g_{1}+g_{2},g_{3})e^{z_{2}D}Y(w,e^{\pi i}z_{2})Y(u,z_{0})v
\nonumber\\
&=&(z_{2}+z_{0})^{k+(g_{1},g_{3})}Y(Y(u,z_{0})v,z_{2})w.
\end{eqnarray}
Then it follows from Proposition \ref{pdlduality}. $\;\;\;\;\Box$

A $V$-{\em module} [DL2] is a vector space 
$W=\coprod_{s\in S}W^{s}$, where $S$ is a $G$-set
equipped with a $\C/2\Z$-valued function $(\cdot,\cdot)$
on $G\times S$ such that
\begin{eqnarray}\label{efunctiongset}
(g_{1}+g_{2},g_{3}+s)=(g_{1},g_{3})+(g_{2},g_{3})+(g_{1},s)+(g_{2},s)
\end{eqnarray}
for $g_{1},g_{2},g_{3}\in G,\;s\in S$, equipped with a
vertex operator map $Y$ from $V$ to $(\End W)\{z\}$ such that the 
axioms (\ref{etruncation}), (\ref{evacuum}) and (\ref{egjacobi0}) 
hold with suitable changes. Sometimes, to distinguish
the vertex operator map $Y$ for a module $W$ from that for
the adjoint module $V$ we use notation $Y_{W}$.

Using the proof of Lemma 2.2 of [DLM]
with a slight modification we get:

\bp{pmodule}
Let $W$ be a $V$-module. Then on $W$,
\begin{eqnarray}
Y(D(v),z)={d\over dz}Y(v,z)\;\;\;\mbox{ for }v\in V.\;\;\;\;\Box
\end{eqnarray}
\ep



At the end of this section we present the following
simple generalization of Lemma 2.3.5 of [Li2]:

\bl{lequivalence}
Let $(V, Y, {\bf 1}, G, c(\cdot,\cdot), (\cdot,\cdot))$ be a
generalized vertex algebra and let $W$ be a $V$-module.
Let $u\in V^{g},\; v\in V^{h},\;n\in \Z,\; u_{(i)}\in V$ for 
$i=1,\dots, k$. If 
\begin{eqnarray}\label{e2.40}
& &(z_{1}-z_{2})^{n+(g,h)}Y_{V}(u,z_{1})Y_{V}(v,z_{2})
-c(g,h)(-1)^{n}(z_{2}-z_{1})^{n+(g,h)}Y_{V}(v,z_{2})Y_{V}(u,z_{1})
\nonumber\\
&=&\sum_{i=0}^{k}Y_{V}(u_{(i)},z_{2})
\left({\partial\over\partial z_{2}}\right)^{i}
\left(z_{1}^{-1}\delta(z_{2}/z_{1})
(z_{2}/z_{1})^{g}\right)
\end{eqnarray}
then 
\begin{eqnarray}
& &(z_{1}-z_{2})^{n+(g,h)}Y_{W}(u,z_{1})Y_{W}(v,z_{2})
-c(g,h)(-1)^{n}(z_{2}-z_{1})^{n+(g,h)}Y_{W}(v,z_{2})Y_{W}(u,z_{1})
\nonumber\\
&=&\sum_{i=0}^{k}Y_{W}(u_{(i)},z_{2})
\left({\partial\over\partial z_{2}}\right)^{i}
\left(z_{1}^{-1}\delta(z_{2}/z_{1})
(z_{2}/z_{1})^{g}\right).
\end{eqnarray}
Furthermore, if $W$ is a faithful module, the converse is also true.
$\;\;\;\;\Box$
\el

\br{rlast}
{\em With Lemma \ref{lequivalence}, we have the following
loose statement: If a generalized vertex algebra $V$ is a module
for a certain ``algebra'' with defining relations of type
(\ref{e2.40}), then any $V$-module is also a module
for this ``algebra.'' Conversely, if $W$ 
is a faithful $V$-module and it is a module for 
a certain ``algebra,'' then $V$
is also a module for the ``algebra.'' Examples for such
``algebras'' are affine Lie algebras, the Virasoro algebra,
affine Griess algebra [FLM] and $Z$-algebras.}
\er

\section{Generalized vertex algebras generated by
parafermion operators}

Throughout this paper, $G$ is an abelian group, 
$(\cdot,\cdot)$ is a symmetric $\Z$-bilinear
$\C/2\Z$-valued form on $G$, $c(\cdot,\cdot)$
is a $\C^{\times}$-valued alternating form on $G$, and
$S$ is a $G$-set equipped with a $\C/2\Z$-valued function 
denoted also by $(\cdot,\cdot)$ on $G\times S$ which satisfies
(\ref{efunctiongset}).

We fix a choice of representatives
$(g,h)$ and $(g,s)$ in $\C$ for $g,h\in G$ and $s\in S$. 
However, it should be observed that the main 
notions and results do not depend on this choice.

Let $W=\coprod_{g\in G}W^{g}$ be a $G$-graded vector space (over $\C$).
By definition,
\begin{eqnarray}
(\End W)\{z\}
=\left\{ a(z)=\sum_{n\in \C}a_{n}z^{-n-1}\;|\; a_{n}\in \End W\right\}.
\end{eqnarray}
Following [DL2], for $g\in G$ we define an operator $z^{g}$ from $W$ to 
$W\{z\}$ by
\begin{eqnarray}
z^{g}\cdot w=z^{(g,s)}w\;\;\;\mbox{ for }w\in W^{s},\; s\in S.
\end{eqnarray}
Note: This operator of course depends on the choice of 
representatives of $(g,s)$.

A formal series $a(z)\in (\End W)\{z\}$ is said to satisfy 
{\em lower truncation condition}
if for every $w\in W$, there exist finitely many complex 
numbers $\a_{1},\dots,\a_{r}$ such that
$$a(z)w\in z^{\a_{1}}W[[z]]+\cdots +z^{\a_{r}}W[[z]].$$
Clearly, all such series form a subspace of $(\End W)\{z\}$.

For $g\in G$, we define $F(W)^{g}$ to be the vector subspace of 
 $(\End W)\{z\}$ consisting of $a(z)$ satisfying the lower truncation condition and 
\begin{eqnarray}
z^{(g,s)}a(z)W^{s}\subset W^{g+s}[[z,z^{-1}]]\;\;\;\mbox{ for }s\in S.
\end{eqnarray}
Note that the notion of $F(W)^{g}$ does not depend on the choice
of representatives $(g,s)$.

Set
\begin{eqnarray}
F(W)=\coprod_{g\in G}F(W)^{g}.
\end{eqnarray}
(It is a direct sum because $W$ is $G$-graded.)
As indicated in the notions of $F(W)^{g}$ and $F(W)$,
$z$ is treated as a dummy variable, i.e., 
$a(z), a(z_{1})$ and $a(z_{2})$ are considered as the same object.

\bd{dglocality}
{\em Formal series $a(z)\in F(W)^{g},\; b(z)\in F(W)^{h}$ 
are said to mutually satisfy
the {\em generalized weak commutativity} if there exists
a nonnegative integer $k$ such that
\begin{eqnarray}\label{eglocality1}
(z_{1}-z_{2})^{k+(g,h)}a(z_{1})b(z_{2})
=(-1)^{k}c(g,h)(z_{2}-z_{1})^{k+(g,h)}b(z_{2})a(z_{1}).
\end{eqnarray}}
\ed
Clearly, this definition is free of the choice of a
representative of $(g,h)$. Note that (\ref{eglocality1}) also holds
if we replace $k$ by any integer greater than $k$.

Recall that $D$ is the linear endomorphism of $(\End W)\{z\}$
such that $D(a(z))=a'(z)$
for $a(z)\in (\End W)\{z\}$, where $a'(z)$ is the formal derivative
of $a(z)$. 
Clearly, $D$ preserves $F(W)$ and its $G$-grading.

\br{rderivative}
{\em  It is easy to see that if $a(z)$ and $b(z)$ mutually satisfy
the generalized weak commutativity, so do $a'(z)$ and $b(z)$.}
\er

A {\em homogeneous parafermion field operator} on $W$ is a
series $a(z)$ of $F(W)^{g}$ for some $g\in G$ that 
satisfies the generalized weak commutativity with itself.

\bd{dvertexproduct}
{\em Let $a(z)\in F(W)^{g},\; b(z)\in F(W)^{h},\;g,h\in G$.
Suppose that $a(z)$ and $b(z)$ satisfy the generalized 
weak commutativity. Then for each $n\in \C$, we define
an element $a(z)_{n}b(z)$ of $(\End W)\{z\}$ by
\begin{eqnarray}\label{e3.6}
(a(z)_{n}b(z))w=\Res_{z_{0}}\Res_{z_{1}}z_{0}^{n}
\left(\frac{z+z_{0}}{z_{1}}\right)^{-(g,s)}X
\end{eqnarray}
for $w\in W^{s},\;s\in S$, where 
\begin{eqnarray}
X&=&z_{0}^{-1}\left(\frac{z_{1}-z}{z_{0}}\right)^{(g,h)}
\delta\left(\frac{z_{1}-z}{z_{0}}\right)a(z_{1})b(z)w\nonumber\\
& &-c(g,h)z_{0}^{-1}
\left(\frac{z-z_{1}}{z_{0}}\right)^{(g,h)}
\delta\left(\frac{z-z_{1}}{-z_{0}}\right)b(z)a(z_{1})w.
\end{eqnarray}}
\ed

\br{rexistence}
{\em Note that because of the generalized weak commutativity,
$z_{0}^{(g,h)}X$, which involves only integral powers of
$z_{0}$, contains only finitely many negative integral powers of $z_{0}$.
Then $a(z)_{n}b(z)$ exists  as a formal series for any
choice of representatives $(g,h), (g,s)$.
Furthermore, since 
\begin{eqnarray}
& &\left(\frac{z_{1}-z}{z_{0}}\right)^{m}
\delta\left(\frac{z_{1}-z}{z_{0}}\right)
=\delta\left(\frac{z_{1}-z}{z_{0}}\right),\\
& &\left(\frac{z-z_{1}}{z_{0}}\right)^{2m}
\delta\left(\frac{z-z_{1}}{-z_{0}}\right)
=\delta\left(\frac{z-z_{1}}{-z_{0}}\right)
\end{eqnarray}
for $m\in \Z$, the expression $X$ does not depend on the choice
of a representative of $(g,h)$. For $m\ge 0$, we also have
\begin{eqnarray*}
& &\left(\frac{z+z_{0}}{z_{1}}\right)^{m}
z_{0}^{-1}\delta\left(\frac{z_{1}-z}{z_{0}}\right)
=z_{0}^{-1}\delta\left(\frac{z_{1}-z}{z_{0}}\right),\\
& &\left(\frac{z+z_{0}}{z_{1}}\right)^{m}
z_{0}^{-1}\delta\left(\frac{z-z_{1}}{-z_{0}}\right)
=z_{0}^{-1}\delta\left(\frac{z-z_{1}}{-z_{0}}\right).
\end{eqnarray*}
Then $(a(z)_{n}b(z))w$ does not depend on the choice
of representatives of $(g,h)$ and $(g,s)$.}
\er

If $n\notin (g,h)+\Z$, the right-hand side of (\ref{e3.6})
does not involve integral powers of $z_{0}$. Therefore,
\begin{eqnarray}
a(z)_{n}b(z)=0\;\;\;\mbox{ for }n\notin (g,h)+\Z. 
\end{eqnarray}
For $n\in (g,h)+\Z$, we have
\begin{eqnarray}\label{eanbzw}
& &(a(z)_{n}b(z))w\nonumber\\
&=&\sum_{i\ge 0}\Res_{z_{1}}{-(g,s)\choose i}z_{1}^{(g,s)}
z^{-(g,s)-i}\cdot \nonumber\\
& &\cdot \left((z_{1}-z)^{n+i}a(z_{1})b(z)w
-(-1)^{n-(g,h)+i}c(g,h)(z-z_{1})^{n+i}
b(z)a(z_{1})w\right).\;\;\;
\end{eqnarray}
Again, because of the generalized weak commutativity,
the sum is really a finite sum. Noticing that 
$a(z_{1})w\in W^{g+s}\{z_{1}\}$ from (\ref{eanbzw}) we have
$$z^{(g,s)+(h,s)}(a(z)_{n}b(z))w\in W^{g+h+s}((z)).$$
Since $(g+h,s)-(g,s)-(h,s)\in 2\Z$, we have
\begin{eqnarray}
z^{(g+h,s)}(a(z)_{n}b(z))w\in W^{g+h+s}((z)).
\end{eqnarray}
This shows that $a(z)_{n}b(z)\in F(W)^{g+h}$. 

To summarize we have:

\bl{lbasic0}
Suppose that $a(z)\in F(W)^{g},\; b(z)\in F(W)^{h}$ 
satisfy the generalized weak commutativity. Then
\begin{eqnarray}
a(z)_{n}b(z)\in F(W)^{g+h}\;\;\;\mbox{ for }n\in \C.
\end{eqnarray}
Furthermore, $a(z)_{n}b(z)=0$ if $n\notin (g,h)+\Z$, and
$a(z)_{n}b(z)=0$ for $n\in (g,h)+\Z$ with $n-(g,h)$ being sufficiently large.$\;\;\;\;\Box$
\el

\br{rdifferent}
{\em Note that in the ordinary untwisted case with $G=0$, or $\Z/2\Z$ 
(cf. [Li2]), 
$a(z)_{n}b(z)$ were well defined 
for all $a(z), b(z)\in (\End W)[[z,z^{-1}]]$.
In the generalized case, 
the definition of $a(z)_{n}b(z)$ requires the generalized
weak commutativity. This is similar to the situation for 
twisted vertex operators in [Li3].}
\er

Write $a(z)_{n}b(z)$ in terms of generating series as 
\begin{eqnarray}
Y(a(z),z_{0})b(z)=\sum_{n\in \C}\left(a(z)_{n}b(z)\right)z_{0}^{-n-1},
\end{eqnarray}
where $z_{0}$ is another formal variable. Then
\begin{eqnarray}\label{eadjoint}
Y(a(z),z_{0})b(z)
={\rm Res}_{z_{1}}\left(\frac{z+z_{0}}{z_{1}}\right)^{-(g,s)}
\cdot X.
\end{eqnarray}
Note that $z_{0}^{(g,h)}Y(a(z),z_{0})b(z)$ involves only integral powers 
of $z_{0}$ and that the powers of $z_{0}$ are truncated from below.
Thus, for $\alpha\in \C$,
$$(z+z_{0})^{\a}Y(a(z),z_{0})b(z)\;\;\;\mbox{exists }$$
in $F(W)\{z_{0},z\}$,
hence
$$\left(\frac{z_{1}-z_{0}}{z}\right)^{\alpha} 
z^{-1}\delta\left(\frac{z_{1}-z_{0}}{z}\right)Y(a(z),z_{0})b(z)
\;\;\;\mbox{exists }$$
in $F(W)\{z_{0},z_{1},z\}$.

As expected we have:

\bp{pjacobirelation}
Let $a(z)\in F(W)^{g},\; b(z)\in F(W)^{h},\;g,h\in G$. 
Suppose that $a(z), b(z)$ mutually 
satisfy the generalized weak commutativity.
Then for $w\in W^{s},\; s\in S$,
\begin{eqnarray}\label{ejacobirelation}
& &\left(\frac{z_{1}-z_{0}}{z}\right)^{-(g,s)} 
z^{-1}\delta\left(\frac{z_{1}-z_{0}}{z}\right)(Y(a(z),z_{0})b(z))w\nonumber\\
&=&z_{0}^{-1}\left(\frac{z_{1}-z}{z_{0}}\right)^{(g,h)}
\delta\left(\frac{z_{1}-z}{z_{0}}\right)a(z_{1})b(z)w\nonumber\\
& &-c(g,h)z_{0}^{-1}
\left(\frac{z-z_{1}}{z_{0}}\right)^{(g,h)}
\delta\left(\frac{z-z_{1}}{-z_{0}}\right)b(z)a(z_{1})w.
\end{eqnarray}
\ep

\pf Let $r$ be a nonnegative integer such that 
$$z_{1}^{r+(g,s)}a(z_{1})w\in W[[z_{1}]],$$
hence
\begin{eqnarray}
\Res_{z_{1}}z_{1}^{r+(g,s)}z_{0}^{-1}
\left(\frac{z-z_{1}}{z_{0}}\right)^{(g,s)}
\delta\left(\frac{z-z_{1}}{-z_{0}}\right)b(z)a(z_{1})w=0.
\end{eqnarray}
Then using (\ref{eadjoint}) and the fundamental properties of 
the delta function, we have
\begin{eqnarray}
& &(z+z_{0})^{r+(g,s)}(Y(a(z),z_{0})b(z))w\nonumber\\
&=&\Res_{z_{1}}(z+z_{0})^{r}z_{1}^{(g,s)}\cdot X\nonumber\\
&=&\Res_{z_{1}}z_{1}^{r+(g,s)}
z_{0}^{-1}\left(\frac{z_{1}-z}{z_{0}}\right)^{(g,h)}
\delta\left(\frac{z_{1}-z}{z_{0}}\right)a(z_{1})b(z)w\nonumber\\
& &-c(g,h) \Res_{z_{1}}z_{1}^{r+(g,s)}
z_{0}^{-1}
\left(\frac{z-z_{1}}{z_{0}}\right)^{(g,h)}
\delta\left(\frac{z-z_{1}}{-z_{0}}\right)b(z)a(z_{1})w\nonumber\\
&=&\Res_{z_{1}}z_{1}^{r+(g,s)}
z_{1}^{-1}\left(\frac{z_{0}+z}{z_{1}}\right)^{-(g,h)}
\delta\left(\frac{z_{0}+z}{z_{1}}\right)a(z_{1})b(z)w\nonumber\\
&=&\Res_{z_{1}}z_{1}^{r+(g,s)-(g,h+s)+(g,h)}
(z_{0}+z)^{-(g,h)}
z_{1}^{-1}\delta\left(\frac{z_{0}+z}{z_{1}}\right)
\left(z_{1}^{(g,h+s)}a(z_{1})b(z)w\right)
\nonumber\\
&=&(z_{0}+z)^{r+(g,s)}a(z_{0}+z)b(z)w,
\end{eqnarray}
noting that $(g,s)-(g,h+s)+(g,h)\in 2\Z$.
Similar to the proof of Proposition \ref{pdlduality}, this generalized 
weak associativity relation
together with the generalized weak commutativity 
relation implies the generalized Jacobi identity. $\;\;\;\;\;\Box$

\bp{pnewlocal}
Let $a(z)\in F(W)^{g_{1}},\; b(z)\in F(W)^{g_{2}},\; 
c(z)\in F(W)^{g_{3}}$ with $g_{1},g_{2},g_{3}\in G$. 
Suppose that $a(z), b(z), c(z)$ mutually 
satisfy the generalized weak commutativity.
Then for $n\in \C$, $a(z)_{n}b(z)$ and $c(z)$ satisfy
the generalized weak commutativity.
\ep

\pf Let $r$ be a positive integer such that the
following identities hold:
\begin{eqnarray*}
& &(z_{1}-z_{2})^{r+(g_{1},g_{3})}a(z_{1})c(z_{2})=(-1)^{r}c(g_{1},g_{3})
(z_{2}-z_{1})^{r+(g_{1},g_{3})}c(z_{2})a(z_{1}),\\
& &(z_{1}-z_{2})^{r+(g_{2},g_{3})}b(z_{1})c(z_{2})
=(-1)^{r}c(g_{2},g_{3})(z_{2}-z_{1})^{r+(g_{2},g_{3})}c(z_{2})b(z_{1}).
\end{eqnarray*}
Let $w\in W^{s},\; s\in S$. 
Using the Jacobi identity relation (\ref{ejacobirelation})
and the above generalized weak commutativity relations we get
\begin{eqnarray}
& &(z_{1}-z_{2})^{r+(g_{1},g_{3})}(z-z_{2})^{r+(g_{2},g_{3})}\cdot \nonumber\\
& &\cdot \left(\frac{z_{1}-z_{0}}{z}\right)^{-(g_{1},g_{3}+s)} 
z^{-1}\delta\left(\frac{z_{1}-z_{0}}{z}\right)
(Y(a(z),z_{0})b(z))c(z_{2})w\nonumber\\
&=&\cdot z_{0}^{-1}\left(\frac{z_{1}-z}{z_{0}}\right)^{(g_{1},g_{2})}
\delta\left(\frac{z_{1}-z}{z_{0}}\right)\cdot \nonumber\\
& &\cdot (z_{1}-z_{2})^{r+(g_{1},g_{3})}(z-z_{2})^{r+(g_{2},g_{3})}
a(z_{1})b(z)c(z_{2})w\nonumber\\
& &-c(g_{1},g_{2})z_{0}^{-1}
\left(\frac{z-z_{1}}{z_{0}}\right)^{(g_{1},g_{2})}
\delta\left(\frac{z-z_{1}}{-z_{0}}\right)\cdot \nonumber\\
& &\cdot (z_{1}-z_{2})^{r+(g_{1},g_{3})}(z-z_{2})^{r+(g_{2},g_{3})}
b(z)a(z_{1})c(z_{2})w\nonumber\\
&=&c(g_{1},g_{3})c(g_{2},g_{3})(z_{2}-z_{1})^{r+(g_{1},g_{3})}
(z_{2}-z)^{r+(g_{2},g_{3})}\cdot \nonumber\\
& &\cdot z_{0}^{-1}\left(\frac{z_{1}-z}{z_{0}}\right)^{(g_{1},g_{2})}
\delta\left(\frac{z_{1}-z}{z_{0}}\right)c(z_{2})a(z_{1})b(z)w\nonumber\\
& &-c(g_{1},g_{3})c(g_{2},g_{3})(z_{2}-z_{1})^{r+(g_{1},g_{3})}
(z_{2}-z)^{r+(g_{2},g_{3})}\cdot \nonumber\\
& &\cdot c(g_{1},g_{2})z_{0}^{-1}
\left(\frac{z-z_{1}}{z_{0}}\right)^{(g_{1},g_{2})}
\delta\left(\frac{z-z_{1}}{-z_{0}}\right)c(z_{2})b(z)a(z_{1})w\nonumber\\
&=&c(g_{1},g_{3})c(g_{2},g_{3})(z_{2}-z_{1})^{r+(g_{1},g_{3})}
(z_{2}-z)^{r+(g_{2},g_{3})}\cdot \nonumber\\
& &\cdot \left(\frac{z_{1}-z_{0}}{z}\right)^{-(g_{1},s)} 
z^{-1}\delta\left(\frac{z_{1}-z_{0}}{z}\right)c(z_{2})
(Y(a(z),z_{0})b(z))w\nonumber\\
&=&c(g_{1}+g_{2},g_{3})(z_{2}-z_{1})^{r+(g_{1},g_{3})}
(z_{2}-z)^{r+(g_{2},g_{3})}
\cdot \nonumber\\
& &\cdot \left(\frac{z_{1}-z_{0}}{z}\right)^{-(g_{1},s)} 
z^{-1}\delta\left(\frac{z_{1}-z_{0}}{z}\right)c(z_{2})(Y(a(z),z_{0})b(z))w.
\end{eqnarray}
Thus
\begin{eqnarray}
& &(z_{1}-z_{2})^{r+(g_{1},g_{3})}(z-z_{2})^{r+(g_{2},g_{3})}\cdot \nonumber\\
& &\cdot \left(\frac{z_{1}-z_{0}}{z}\right)^{-(g_{1},g_{3})} 
z^{-1}\delta\left(\frac{z_{1}-z_{0}}{z}\right)(Y(a(z),z_{0})b(z))c(z_{2})w
\nonumber\\
&=&c(g_{1}+g_{2},g_{3})(z_{2}-z_{1})^{r+(g_{1},g_{3})}
(z_{2}-z)^{r+(g_{2},g_{3})}
\cdot \nonumber\\
& &\cdot z^{-1}\delta\left(\frac{z_{1}-z_{0}}{z}\right)
c(z_{2})(Y(a(z),z_{0})b(z))w,
\end{eqnarray}
noting that $(g_{1},g_{3}+s)-(g_{1},g_{3})-(g_{1},s)\in 2\Z$.
Using the fundamental properties of delta-function  we get
\begin{eqnarray}\label{eleft-left}
& &(z+z_{0}-z_{2})^{r+(g_{1},g_{3})}(z-z_{2})^{r+(g_{2},g_{3})}
\cdot \nonumber\\
& &\cdot 
z_{1}^{-1}\delta\left(\frac{z+z_{0}}{z_{1}}\right)(Y(a(z),z_{0})b(z))c(z_{2})w
\nonumber\\
&=&c(g_{1}+g_{2},g_{3})(z_{2}-z-z_{0})^{r+(g_{1},g_{3})}
(z_{2}-z)^{r+(g_{2},g_{3})}\cdot \nonumber\\
& &\cdot 
z_{1}^{-1}\delta\left(\frac{z+z_{0}}{z_{1}}\right)c(z_{2})(Y(a(z),z_{0})b(z))w.
\end{eqnarray}
Taking $\Res_{z_{1}}$ from (\ref{eleft-left}) we obtain
\begin{eqnarray}\label{eleft-leftnew}
& &(z+z_{0}-z_{2})^{r+(g_{1},g_{3})}(z-z_{2})^{r+(g_{2},g_{3})}
(Y(a(z),z_{0})b(z))c(z_{2})w\nonumber\\
&=&c(g_{1}+g_{2},g_{3})(z_{2}-z-z_{0})^{r+(g_{1},g_{3})}
(z_{2}-z)^{r+(g_{2},g_{3})}
c(z_{2})(Y(a(z),z_{0})b(z))w.\;\;
\end{eqnarray}
Let $n\in \C$ be arbitrarily fixed. 
Since $a(z)_{n}b(z)=0$ for $n\notin (g_{1},g_{2})+\Z$,
we only need to consider $n\in (g_{1},g_{2})+\Z$.
Let $N$ be a fixed nonnegative integer such that
$a(z)_{m}b(z)=0$ for $m\ge N+n$, so that
\begin{eqnarray}
\Res_{z_{0}}z_{0}^{N+n+i}Y(a(z),z_{0})b(z)=0
\end{eqnarray}
for $i\in \N$. We may replace $r$ by $r+N$ so that we may assume 
that $r\ge N$ and that $r+n-(g_{1},g_{2})>0$.
Then using (\ref{eleft-leftnew}) we obtain
\begin{eqnarray}
& &(z-z_{2})^{3r+(g_{1},g_{3})+(g_{2},g_{3})}(a(z)_{n}b(z))c(z_{2})w\nonumber\\
&=&\Res_{z_{0}}z_{0}^{n}(z-z_{2})^{3r+(g_{1},g_{3})+(g_{2},g_{3})}
(Y(a(z),z_{0})b(z))c(z_{2})w\nonumber\\
&=&\Res_{z_{0}}\sum_{i\ge 0}(-1)^{i}{2r+(g_{1},g_{3})\choose i}z_{0}^{n+i}
(z-z_{2}+z_{0})^{2r+(g_{1},g_{3})-i}\cdot \nonumber\\
& &\cdot(z-z_{2})^{r+(g_{2},g_{3})}(Y(a(z),z_{0})b(z))c(z_{2})w\nonumber\\
&=&\Res_{z_{0}}\sum_{i=0}^{N}(-1)^{i}{2r+(g_{1},g_{3})\choose i}z_{0}^{n+i}
\cdot \nonumber\\
& &\cdot (z-z_{2}+z_{0})^{2r+(g_{1},g_{3})-i}
(z-z_{2})^{r+(g_{2},g_{3})}(Y(a(z),z_{0})b(z))c(z_{2})w\nonumber\\
&=&
\Res_{z_{0}}\sum_{i=0}^{N}(-1)^{i}{2r+(g_{1},g_{3})\choose i}z_{0}^{n+i}
c(g_{1}+g_{2},g_{3})(-1)^{r-i}\cdot\nonumber\\
& &\cdot(z_{2}-z-z_{0})^{2r+(g_{1},g_{3})-i}(z+z_{0})^{(g_{1},g_{3})}
(z_{2}-z)^{r+(g_{2},g_{3})}c(z_{2})(Y(a(z),z_{0})b(z))w\nonumber\\
&=&
\Res_{z_{0}}\sum_{i\ge 0}(-1)^{r}{2r+(g_{1},g_{3})\choose i}z_{0}^{n+i}
 c(g_{1}+g_{2},g_{3})\cdot\nonumber\\
& &\cdot(z_{2}-z-z_{0})^{2r+(g_{1},g_{3})-i}(z+z_{0})^{(g_{1},g_{3})}
(z_{2}-z)^{r+(g_{2},g_{3})}c(z_{2})(Y(a(z),z_{0})b(z))w\nonumber\\
&=&(-1)^{r}\Res_{z_{0}}z_{0}^{n}c(g_{1}+g_{2},g_{3})
(z_{2}-z)^{2r+(g_{1},g_{3})+(g_{2},g_{3})}
c(z_{2})(Y(a(z),z_{0})b(z))w\nonumber\\
&=&(-1)^{r}c(g_{1}+g_{2},g_{3})(z_{2}-z)^{3r+(g_{1},g_{3})+(g_{2},g_{3})}
c(z_{2})(a(z)_{n}b(z))w.
\end{eqnarray}
Since $(g_{1}+g_{2},g_{3})-(g_{1},g_{3})-(g_{2},g_{3})\in 2\Z$,
there exist $k,k'\in \N$ such that
$$2k+(g_{1},g_{3})+(g_{2},g_{3})=(g_{1}+g_{2},g_{3})+2k'.$$
Then
\begin{eqnarray}
& &(z-z_{2})^{3r+2k'+(g_{1}+g_{2},g_{3})}(a(z)_{n}b(z))c(z_{2})w
\nonumber\\
&=&(-1)^{r}c(g_{1}+g_{2},g_{3})(z_{2}-z)^{3r+2k'+(g_{1},+g_{2},g_{3})}
c(z_{2})(a(z)_{n}b(z))w.
\end{eqnarray}
This proves that $a(z)_{n}b(z)$ and $c(z)$ mutually satisfy 
the generalized weak commutativity. $\;\;\;\;\Box$

\bd{dpartialalgebra}
{\em A $G$-graded subspace $A$ of $F(W)$ is called 
a {\em generalized vertex pre-algebra} if every pair of 
homogeneous elements of $A$
satisfy the generalized weak commutativity.} 
\ed

For homogeneous $a(z), b(z)\in A$, $a(z)_{n}b(z)$ was defined for $n\in \C$.
Using linearity, we define $a(z)_{n}b(z)$ for all $a(z),b(z)\in A$, so that
$Y(a(z),z_{0})b(z)$ is defined for all $a(z),b(z)\in A$.
Furthermore, $A$ is 
said to be {\em closed} if $a(z)_{n}b(z)\in A$ for $a(z),b(z)\in A, n\in \C$.
Since the identity operator $I(z)=\mbox{id}_{W}$ (independent of $z$) and 
any element of $F(W)$ mutually satisfy the generalized weak 
commutativity, any maximal generalized vertex pre-algebra
contains $I(z)$. 
In view of Proposition \ref{pnewlocal} we immediately have:

\bc{cmaximal}
Any maximal generalized vertex pre-algebra $V$ contains the identity 
operator $I(z)$ and it is closed and $D$-stable.$\;\;\;\;\Box$
\ec

For the rest of this section,
$V$ will be a fixed closed generalized vertex pre-algebra containing 
the identity operator $I(z)$ on $W$. The same proof of Lemma 3.10 of [Li3]
gives the following results:

\bl{lbasic1}
For $a(z)\in V$, we have
\begin{eqnarray}
& &Y(I(z),z_{0})a(z)=a(z),\\
& &\;Y(a(z),z_{0})I(z)=e^{z_{0}{\partial\over\partial z}}a(z)
=a(z+z_{0}). \;\;\;\;\Box
\end{eqnarray}
\el


Furthermore, we have:

\bl{ld-derivative}
Suppose that $a(z)\in F(W)^{g},\; b(z)\in F(W)^{h}$ 
mutually satisfy the generalized weak commutativity.
Then
\begin{eqnarray}
& &{\partial\over \partial z_{0}}Y(a(z),z_{0})b(z)=Y(a'(z),z_{0})b(z)=
Y(D(a(z)),z_{0})b(z),\\
& &[D, Y(a(z),z_{0})]b(z)={\partial\over\partial z_{0}}Y(a(z),z_{0})b(z).
\end{eqnarray}
\el

{\bf Proof}. Let $w\in W^{s},\; s\in S$.
Recall that the generalized Jacobi relation (\ref{ejacobirelation}) holds.
We also have
\begin{eqnarray}\label{ea'jacobirelation}
& &\left(\frac{z_{1}-z_{0}}{z}\right)^{-(g,s)} 
z^{-1}\delta\left(\frac{z_{1}-z_{0}}{z}\right)(Y(a'(z),z_{0})b(z))w\nonumber\\
&=&z_{0}^{-1}\left(\frac{z_{1}-z}{z_{0}}\right)^{(g,h)}
\delta\left(\frac{z_{1}-z}{z_{0}}\right)a'(z_{1})b(z)w\nonumber\\
& &-c(g,h)z_{0}^{-1}
\left(\frac{z-z_{1}}{z_{0}}\right)^{(g,h)}
\delta\left(\frac{z-z_{1}}{-z_{0}}\right)b(z)a'(z_{1})w.
\end{eqnarray}
Let $X_{L}$ and $X_{R}$ be the term on the left-hand side and 
the term on the right-hand side of (\ref{ejacobirelation}).
Using (\ref{ea'jacobirelation}) we get
\begin{eqnarray}
& &\left(\frac{z+z_{0}}{z_{1}}\right)^{(g,s)} 
z_{1}^{-1}\delta\left(\frac{z+z_{0}}{z_{1}}\right)
(Y(a'(z),z_{0})b(z))w\nonumber\\
&=&\left(\frac{z_{1}-z_{0}}{z}\right)^{-(g,s)} 
z^{-1}\delta\left(\frac{z_{1}-z_{0}}{z}\right)(Y(a'(z),z_{0})b(z))w\nonumber\\
&=&{\partial\over\partial z_{1}}X_{R}
-a(z_{1})b(z)w{\partial\over \partial z_{1}}
\left(z_{0}^{-1}\left(\frac{z_{1}-z}{z_{0}}\right)^{(g,h)}
\delta\left(\frac{z_{1}-z}{z_{0}}\right)\right)\nonumber\\
& &+c(g,h)b(z)a(z_{1})w
{\partial\over \partial z_{1}}
\left(z_{0}^{-1}\left(\frac{z-z_{1}}{z_{0}}\right)^{(g,h)}
\delta\left(\frac{z-z_{1}}{-z_{0}}\right)\right)\nonumber\\
&=&{\partial\over\partial z_{1}}X_{L}
+a(z_{1})b(z)w{\partial\over \partial z_{0}}
\left(z_{0}^{-1}\left(\frac{z_{1}-z}{z_{0}}\right)^{(g,h)}
\delta\left(\frac{z_{1}-z}{z_{0}}\right)\right)\nonumber\\
& &-c(g,h)b(z)a(z_{1})w
{\partial\over \partial z_{0}}
\left(z_{0}^{-1}\left(\frac{z-z_{1}}{z_{0}}\right)^{(g,h)}
\delta\left(\frac{z-z_{1}}{-z_{0}}\right)\right)\nonumber\\
&=&{\partial\over\partial z_{1}}X_{L}+{\partial\over\partial z_{0}}X_{R}
\nonumber\\
&=&{\partial\over\partial z_{1}}X_{L}+{\partial\over\partial z_{0}}X_{L}.
\end{eqnarray}
We are using the fact
\begin{eqnarray}\label{epartial}
{\partial\over\partial z_{1}}
\left(\left(\frac{z-z_{1}}{z_{0}}\right)^{\a} 
z_{0}^{-1}\delta\left(\frac{z-z_{1}}{z_{0}}\right)\right)
=-{\partial\over\partial z_{0}}
\left(\left(\frac{z-z_{1}}{z_{0}}\right)^{\a} 
z_{0}^{-1}\delta\left(\frac{z-z_{1}}{z_{0}}\right)\right).
\end{eqnarray}
for $\a\in \C$.
Multiplying by $z_{1}^{(g,s)}$, then taking $\Res_{z_{1}}$ 
and using a variant of (\ref{epartial}) we get
\begin{eqnarray}
& &(z+z_{0})^{(g,s)}(Y(a'(z),z_{0})b(z))w\nonumber\\
&=&\Res_{z_{1}}z_{1}^{(g,s)}{\partial\over\partial z_{1}}X_{L}
+\Res_{z_{1}}z_{1}^{(g,s)}{\partial\over\partial z_{0}}X_{L}
\nonumber\\
&=&\Res_{z_{1}}(Y(a(z),z_{0})b(z))w z_{1}^{(g,s)}
{\partial\over\partial z_{1}}
\left(\left(\frac{z+z_{0}}{z_{1}}\right)^{(g,s)} 
z_{1}^{-1}\delta\left(\frac{z+z_{0}}{z_{1}}\right)\right)\nonumber\\
& &+\Res_{z_{1}}(Y(a(z),z_{0})b(z))w z_{1}^{(g,s)}
{\partial\over\partial z_{0}}
\left(\left(\frac{z+z_{0}}{z_{1}}\right)^{(g,s)} 
z_{1}^{-1}\delta\left(\frac{z+z_{0}}{z_{1}}\right)\right)\nonumber\\
& &+\Res_{z_{1}}z_{1}^{(g,s)}
\left(\frac{z+z_{0}}{z_{1}}\right)^{(g,s)} 
z_{1}^{-1}\delta\left(\frac{z+z_{0}}{z_{1}}\right)
{\partial\over\partial z_{0}}(Y(a(z),z_{0})b(z))w\nonumber\\
&=&\Res_{z_{1}}z_{1}^{(g,s)}
\left(\frac{z+z_{0}}{z_{1}}\right)^{(g,s)} 
z_{1}^{-1}\delta\left(\frac{z+z_{0}}{z_{1}}\right)
{\partial\over\partial z_{0}}(Y(a(z),z_{0})b(z))w\nonumber\\
&=&(z+z_{0})^{(g,s)}{\partial\over\partial z_{0}}(Y(a(z),z_{0})b(z))w.
\end{eqnarray}
Multiplying by $(z+z_{0})^{-(g,s)}$ from left we get the first identity.

The second identity follows from a similar argument
and the fact:
\begin{eqnarray}
{\partial\over\partial z_{1}}
\left(\left(\frac{z_{1}-z}{z_{0}}\right)^{\a} 
z_{0}^{-1}\delta\left(\frac{z_{1}-z}{z_{0}}\right)\right)
=-{\partial\over\partial z}
\left(\left(\frac{z_{1}-z}{z_{0}}\right)^{\a} 
z_{0}^{-1}\delta\left(\frac{z_{1}-z}{z_{0}}\right)\right)
\end{eqnarray}
for $\a\in \C$.$\;\;\;\;\Box$

\bp{padjointlocal}
Let $V$ be a closed generalized 
vertex pre-algebra of parafermion operators on $W$. 
Let $a(z)\in V^{g_{1}},\;b(z)\in V^{g_{2}}$ and let $r\in \N$ 
be such that
\begin{eqnarray}
(z_{1}-z_{2})^{r+(g_{1},g_{2})}a(z_{1})b(z_{2})
=(-1)^{r}c(g_{1},g_{2})(z_{2}-z_{1})^{r+(g_{1},g_{2})}b(z_{2})a(z_{1}).
\end{eqnarray}
Then
\begin{eqnarray}\label{eadkointgwc}
& &(x_{1}-x_{2})^{r+(g_{1},g_{2})}Y(a(z),x_{1})Y(b(z),x_{2})\nonumber\\
&=&(-1)^{r}c(g_{1},g_{2})(x_{2}-x_{1})^{r+(g_{1},g_{2})}
Y(b(z),x_{2})Y(a(z),x_{1}),
\end{eqnarray}
acting on $V$.
\ep

{\bf Proof}. Let $c(z)\in V^{g_{3}},\; w\in W^{s},\; g_{3}\in G,\; s\in S$. 
Using the generalized Jacobi relation (\ref{ejacobirelation})
and the fundamental properties of the delta-function we get
\begin{eqnarray}\label{emain1}
& &\left(\frac{z+x_{1}}{z_{1}}\right)^{(g_{1},s)}
z_{1}^{-1}\delta\left(\frac{z+x_{1}}{z_{1}}\right)
\left(\frac{z+x_{2}}{z_{2}}\right)^{(g_{2},s)}
z_{2}^{-1}\delta\left(\frac{z+x_{2}}{z_{2}}\right)\cdot \nonumber\\
& &\cdot (Y(a(z),x_{1})Y(b(z),x_{2})c(z))w
\nonumber\\
&=&\left(\frac{z_{1}-x_{1}}{z}\right)^{-(g_{1},s)}
z^{-1}\delta\left(\frac{z_{1}-x_{1}}{z}\right)
\left(\frac{z_{2}-x_{2}}{z}\right)^{-(g_{2},s)}
z^{-1}\delta\left(\frac{z_{2}-x_{2}}{z}\right)\cdot \nonumber\\
& &\cdot (Y(a(z),x_{1})Y(b(z),x_{2})c(z))w
\nonumber\\
&=&\left(\frac{z_{2}-x_{2}}{z}\right)^{-(g_{2},s)}
z^{-1}\delta\left(\frac{z_{2}-x_{2}}{z}\right)\cdot \nonumber\\
& &\cdot x_{1}^{-1}\left(\frac{z_{1}-z}{x_{1}}\right)^{(g_{1},g_{2}+g_{3})}
\delta\left(\frac{z_{1}-z}{x_{1}}\right)
a(z_{1})(Y(b(z),x_{2})c(z))w\nonumber\\
& &-c(g_{1},g_{2}+g_{3})
\left(\frac{z_{2}-x_{2}}{z}\right)^{-(g_{2},s)}
z^{-1}\delta\left(\frac{z_{2}-x_{2}}{z}\right)\cdot \nonumber\\
& &\cdot x_{1}^{-1}\left(\frac{z-z_{1}}{x_{1}}\right)^{(g_{1},g_{2}+g_{3})}
\delta\left(\frac{z-z_{1}}{-x_{1}}\right)(Y(b(z),x_{2})c(z))a(z_{1})w
\nonumber\\
&=&x_{1}^{-1}\left(\frac{z_{1}-z}{x_{1}}\right)^{(g_{1},g_{2}+g_{3})}
\delta\left(\frac{z_{1}-z}{x_{1}}\right)\cdot \nonumber\\
& &\cdot x_{2}^{-1}\left(\frac{z_{2}-z}{x_{2}}\right)^{(g_{2},g_{3})}
\delta\left(\frac{z_{2}-z}{x_{2}}\right)a(z_{1})b(z_{2})c(z)w\nonumber\\
& &-x_{1}^{-1}\left(\frac{z_{1}-z}{x_{1}}\right)^{(g_{1},g_{2}+g_{3})}
\delta\left(\frac{z_{1}-z}{x_{1}}\right)\cdot \nonumber\\
& &\cdot c(g_{2},g_{3})
x_{2}^{-1}\left(\frac{z-z_{2}}{x_{2}}\right)^{(g_{2},g_{3})}
\delta\left(\frac{z-z_{2}}{-x_{2}}\right)a(z_{1})c(z)b(z_{2})w \nonumber\\
& &-c(g_{1},g_{2}+g_{3})
x_{1}^{-1}\left(\frac{z-z_{1}}{x_{1}}\right)^{(g_{1},g_{2}+g_{3})}
\delta\left(\frac{z-z_{1}}{-x_{1}}\right)\cdot \nonumber\\
& &\cdot \left(\frac{z_{2}-x_{2}}{z}\right)^{-(g_{2},s)}
z^{-1}\delta\left(\frac{z_{2}-x_{2}}{z}\right)
(Y(b(z),x_{2})c(z))a(z_{1})w.
\end{eqnarray}
Let $p,q\in \N$ be such that
\begin{eqnarray}\label{estcondition}
z_{1}^{p+(g_{1},s)}a(z_{1})w\in W[[z_{1}]],\;\;\;
z_{2}^{q+(g_{2},s)}b(z_{2})w\in W[[z_{2}]].
\end{eqnarray}
Notice that for $j\in \N$ we have
\begin{eqnarray}\label{eproperty2}
(z+x_{3})^{j}A=z_{3}^{j}A,
\end{eqnarray}
where $A$ is one of the three delta-functions 
$x_{3}^{-1}\delta\left(\frac{z_{3}-z}{x_{3}}\right),\; 
x_{3}^{-1}\delta\left(\frac{z-z_{3}}{-x_{3}}\right)$ and 
$z^{-1}\delta\left(\frac{z_{3}-x_{3}}{z}\right)$. 
Applying 
$\Res_{z_{1}}\Res_{z_{2}}(z+x_{1})^{s}(z+x_{2})^{q}
z_{1}^{(g_{1},s)}z_{2}^{(g_{2},s)}$ to (\ref{emain1}), then using
(\ref{eproperty2}) and (\ref{estcondition}) we get
\begin{eqnarray}
& &(z+x_{1})^{p+(g_{1},h)}(z+x_{2})^{q+(g_{2},s)}
(Y(a(z),x_{1})Y(b(z),x_{2})c(z))w\nonumber\\
&=&\Res_{z_{1}}\Res_{z_{2}}z_{1}^{p+(g_{1},h)}z_{2}^{q+(g_{2},s)}
x_{1}^{-1}\left(\frac{z_{1}-z}{x_{1}}\right)^{(g_{1},g_{2}+g_{3})}
\delta\left(\frac{z_{1}-z}{x_{1}}\right)\cdot \nonumber\\
& &\cdot x_{2}^{-1}\left(\frac{z_{2}-z}{x_{2}}\right)^{(g_{2},g_{3})}
\delta\left(\frac{z_{2}-z}{x_{2}}\right)a(z_{1})b(z_{2})c(z)w\nonumber\\
& &-\Res_{z_{1}}\Res_{z_{2}}z_{1}^{p+(g_{1},s)}
x_{1}^{-1}\left(\frac{z_{1}-z}{x_{1}}\right)^{(g_{1},g_{2}+g_{3})}
\delta\left(\frac{z_{1}-z}{x_{1}}\right)\cdot \nonumber\\
& &\cdot c(g_{2},g_{3})
x_{2}^{-1}\left(\frac{z-z_{2}}{x_{2}}\right)^{(g_{2},g_{3})}
\delta\left(\frac{z-z_{2}}{-x_{2}}\right)a(z_{1})c(z)
\left(z_{2}^{q+(g_{2},s)}b(z_{2})w\right) \nonumber\\
& &-\Res_{z_{1}}\Res_{z_{2}}z_{2}^{q+(g_{2},s)}
c(g_{1},g_{2}+g_{3})
x_{1}^{-1}\left(\frac{z-z_{1}}{x_{1}}\right)^{(g_{1},g_{2}+g_{3})}
\delta\left(\frac{z-z_{1}}{-x_{1}}\right)\cdot \nonumber\\
& &\cdot \left(\frac{z_{2}-x_{2}}{z}\right)^{-(g_{2},s)}
z^{-1}\delta\left(\frac{z_{2}-x_{2}}{z}\right)
(Y(b(z),x_{2})c(z))\left(z_{1}^{p+(g_{1},s)}a(z_{1})w\right)\nonumber\\
&=&\Res_{z_{1}}\Res_{z_{2}}z_{1}^{p+(g_{1},s)}z_{2}^{q+(g_{2},s)}
x_{1}^{-1}\left(\frac{z_{1}-z}{x_{1}}\right)^{(g_{1},g_{2}+g_{3})}
\delta\left(\frac{z_{1}-z}{x_{1}}\right)\cdot \nonumber\\
& &\cdot x_{2}^{-1}\left(\frac{z_{2}-z}{x_{2}}\right)^{(g_{2},g_{3})}
\delta\left(\frac{z_{2}-z}{x_{2}}\right)a(z_{1})b(z_{2})c(z)w.
\end{eqnarray}
Notice that
\begin{eqnarray}
& &(x_{1}-x_{2})^{r+(g_{1},g_{2})}
x_{1}^{-1}\delta\left(\frac{z_{1}-z}{x_{1}}\right)
x_{2}^{-1}\delta\left(\frac{z_{2}-z}{x_{2}}\right)\nonumber\\
&=&(x_{1}-z_{2}+z)^{r+(g_{1},g_{2})}
x_{1}^{-1}\delta\left(\frac{z_{1}-z}{x_{1}}\right)
x_{2}^{-1}\delta\left(\frac{z_{2}-z}{x_{2}}\right)\nonumber\\
&=&\left(\frac{z_{1}-z}{x_{1}}\right)^{-r-(g_{1},g_{2})}
(z_{1}-z_{2})^{r+(g_{1},g_{2})}
x_{1}^{-1}\delta\left(\frac{z_{1}-z}{x_{1}}\right)
x_{2}^{-1}\delta\left(\frac{z_{2}-z}{x_{2}}\right)\nonumber\\
&=&\left(\frac{z_{1}-z}{x_{1}}\right)^{-(g_{1},g_{2})}
(z_{1}-z_{2})^{r+(g_{1},g_{2})}
x_{1}^{-1}\delta\left(\frac{z_{1}-z}{x_{1}}\right)
x_{2}^{-1}\delta\left(\frac{z_{2}-z}{x_{2}}\right).
\end{eqnarray}
Then
\begin{eqnarray}
& &(z+x_{1})^{p+(g_{1},s)}(z+x_{2})^{q+(g_{2},s)}
(x_{1}-x_{2})^{r+(g_{1},g_{2})}(Y(a(z),x_{1})Y(b(z),x_{2})c(z))w\nonumber\\
&=&\Res_{z_{1}}\Res_{z_{2}}z_{1}^{p+(g_{1},s)}z_{2}^{q+(g_{2},s)}
x_{1}^{-1}\left(\frac{z_{1}-z}{x_{1}}\right)^{(g_{1},g_{3})}
\delta\left(\frac{z_{1}-z}{x_{1}}\right)\cdot \nonumber\\
& &\cdot x_{2}^{-1}\left(\frac{z_{2}-z}{x_{2}}\right)^{(g_{2},g_{3})}
\delta\left(\frac{z_{2}-z}{x_{2}}\right)
(z_{1}-z_{2})^{r+(g_{1},g_{2})}a(z_{1})b(z_{2})c(z)w.
\end{eqnarray}
Using the obvious symmetry, we have
\begin{eqnarray}
& &(z+x_{1})^{p+(g_{1},s)}(z+x_{2})^{q+(g_{2},s)}
(x_{2}-x_{1})^{r+(g_{1},g_{2})}(Y(b(z),x_{2})Y(a(z),x_{1})c(z))w\nonumber\\
&=&\Res_{z_{1}}\Res_{z_{2}}z_{1}^{p+(g_{1},s)}z_{2}^{q+(g_{2},s)}
x_{2}^{-1}\left(\frac{z_{2}-z}{x_{2}}\right)^{(g_{2},g_{3})}
\delta\left(\frac{z_{2}-z}{x_{2}}\right)\cdot \nonumber\\
& &\cdot x_{1}^{-1}\left(\frac{z_{1}-z}{x_{1}}\right)^{(g_{1},g_{3})}
\delta\left(\frac{z_{1}-z}{x_{1}}\right)
(z_{2}-z_{1})^{r+(g_{1},g_{2})}b(z_{2})a(z_{1})c(z)w.
\end{eqnarray}
Therefore
\begin{eqnarray}
& &(z+x_{1})^{p+(g_{1},s)}(z+x_{2})^{q+(g_{2},s)}
(x_{1}-x_{2})^{r+(g_{1},g_{2})}(Y(a(z),x_{1})Y(b(z),x_{2})c(z))w\nonumber\\
&=&(-1)^{r}c(g_{1},g_{2})(z+x_{1})^{p+(g_{1},s)}(z+x_{2})^{q+(g_{2},s)}
(x_{2}-x_{1})^{r+(g_{1},g_{2})}\cdot\nonumber\\
& &\cdot (Y(b(z),x_{2})Y(a(z),x_{1})c(z))w.
\end{eqnarray}
Multiplying both sides by $(z+x_{1})^{-p-(g_{1},s)}(z+x_{2})^{-q-(g_{2},s)}$
we obtain
\begin{eqnarray}
& &(x_{1}-x_{2})^{r+(g_{1},g_{2})}(Y(a(z),x_{1})Y(b(z),x_{2})c(z))w\nonumber\\
&=&(-1)^{r}c(g_{1},g_{2})
(x_{2}-x_{1})^{r+(g_{1},g_{2})}(Y(b(z),x_{2})Y(a(z),x_{1})c(z))w.
\end{eqnarray}
Then the generalized weak commutativity (\ref{eadkointgwc}) follows
immediately.$\;\;\;\;\Box$

Now, we are in a position to present our main theorem:

\bt{tclosedalgebra}
Let $V$ be a closed generalized vertex pre-algebra of parafermionic 
field operators on $W$, containing the identity operator $I(z)$.
Then $V$ is a generalized vertex algebra 
and $W$ is a canonical $V$-module with $Y(a(z),z_{1})=a(z_{1})$ 
for $a(z)\in V$.
\et

\pf It follows from Propositions \ref{pdlgvoa}, \ref{padjointlocal} and 
Lemmas \ref{lbasic0}, \ref{lbasic1} and \ref{ld-derivative} that $(V,Y)$ is a
generalized vertex algebra. It follows from Proposition \ref{pjacobirelation}
that $W$ is a $V$-module under the natural action.
$\;\;\;\;\Box$

The following is a very useful consequence:

\bt{talgebra}
Let $\Gamma$ be a set of homogeneous parafermion field operators on $W$ 
that satisfy the generalized weak commutativity. Then 
the subspace $\<\Gamma\>$ of $F(W)$, linearly spanned by
\begin{eqnarray}\label{eelementslist}
a^{(1)}(z)_{n_{1}}\cdots a^{(r)}(z)_{n_{r}}I(z)
\end{eqnarray}
for $r\in \N,\; a^{(i)}(z)\in \Gamma,\; n_{1},\dots, n_{r}\in \C$,
equipped with the vertex operator map $Y$ is a 
generalized vertex algebra with $W$ as a natural module.
\et

\pf By Zorn's lemma, there exists a maximal 
generalized vertex pre-algebra $V$ containing $\Gamma$. 
By Corollary \ref{cmaximal}, $V$ is closed and contains $I(z)$.
By Theorem \ref{tclosedalgebra}, $V$ is a generalized vertex algebra.
It follows from Proposition 14.8 of [DL2] that $\<\Gamma\>$ 
is a generalized vertex subalgebra of $V$. Clearly,
$\<\Gamma\>$ does not depend on the choice of $V$
and it is the intersection of all closed generalized 
vertex pre-algebra containing $\Gamma$ and the identity operator 
$I(z)$. $\;\;\;\;\Box$

\bl{lvacuumlikevector}
Let $V$ be a generalized vertex algebra and let $U$ be a
graded subspace which generates $V$. Let $W$ be a 
$V$-module and let $e\in W^{0}$ such that
$Y(u,z)e\in V[[z]]$ for $u\in U$. Then
\begin{eqnarray}
Y(v,z)e\in W[[z]]\;\;\;\mbox{ for all }v\in V.
\end{eqnarray}
\el

\pf Let $V'$ be the collection of all $v\in V$ such that $Y(v,z)e\in V[[z]]$.
Clearly, $V'$ is a graded subspace of $V$ containing ${\bf 1}$ and $U$.
Let $u\in (V')^{g},\; v\in (V')^{h},\; g,h\in G$. Then
$Y(u,z)e, Y(v,z)e\in V[[z]]$.
Applying the generalized Jacobi identity to $e$, then taking $\Res_{z_{1}}$
 we get
\begin{eqnarray}
& &Y(Y(u,z_{0})v,z_{2})e\nonumber\\
&=&\Res_{z_{1}}z_{0}^{-1}\left(\frac{z_{1}-z_{2}}{z_{0}}\right)^{(g,h)}
\delta\left(\frac{z_{1}-z_{2}}{z_{0}}\right)Y(u,z_{1})Y(v,z_{2})e\nonumber\\
& &- \Res_{z_{1}}c(g,h)z_{0}^{-1}\left(\frac{z_{2}-z_{1}}{z_{0}}\right)^{(g,h)}
\delta\left(\frac{z_{2}-z_{1}}{-z_{0}}\right)Y(v,z_{2})Y(u,z_{1})e\nonumber\\
&=&\Res_{z_{1}}z_{0}^{-1}\left(\frac{z_{1}-z_{2}}{z_{0}}\right)^{(g,h)}
\delta\left(\frac{z_{1}-z_{2}}{z_{0}}\right)Y(u,z_{1})Y(v,z_{2})e\nonumber\\
&\in& (V((z_{0})))[[z_{2}]].
\end{eqnarray}
Then
\begin{eqnarray}
Y(u_{n}v,z_{2})e\in V[[z_{2}]]\;\;\;\mbox{ for }n\in \C.
\end{eqnarray}
Thus $u_{n}v\in V'$. Thus, $V'$ is a subalgebra of $V$, 
containing ${\bf 1}$ and $U$. Consequently, $V'=V$.
The proof is complete.$\;\;\;\;\Box$

The same proof of Proposition 3.4 of [Li1] gives:

\bp{pvacuumlikevector2}
Let $V$ be a generalized vertex algebra, let $W$ be a 
$V$-module and let $e\in W^{0}$ such that
$Y(v,z)e\in V[[z]]$ for $v\in V$. Then 
\begin{eqnarray}
Y(v,z)e=e^{zL(-1)}v_{-1}e \;\;\;\mbox{ for }v\in V,
\end{eqnarray}
and the linear map
\begin{eqnarray}
f: & &V\rightarrow W\nonumber\\
& &v\mapsto v_{-1}e
\end{eqnarray}
is a $V$-homomorphism. Furthermore, if $W$ is faithful and generated by $e$, then
$f$ is an isomorphism.$\;\;\;\;\Box$
\ep

Now we have the following result the first part of which is a generalization 
of a result obtained in [FKRW] and [PM]:

\bt{tfkrw}
Let $V$ be a $G$-graded vector space and ${\bf 1}\in V^{0}$
and let $U$ a generalized vertex pre-algebra of parafermions on $V$
such that $\psi(z){\bf 1}\in V[[z]]$ and such that $V$ is generated from
${\bf 1}$ by all component operators of $\psi(z)z^{g}$ 
for $\psi(z)\in U^{g},\; g\in G$.
Then there exists a unique generalized vertex algebra structure on $V$ such
that ${\bf 1}$ is the vacuum vector and that $Y(a,z)=a(z)$ for $a\in U$.
Furthermore, this generalized vertex algebra $V$ is isomorphic to
the generalized vertex algebra generated by $U$ inside $(\End V)\{z\}$.
\et

\pf Let $\bar{U}$ be the generalized vertex algebra generated by $U$
inside $F(V)$. Then $V$ is a $\bar{U}$-module and $Y(u,z){\bf 1}\in V[[z]]$ 
for $u\in U$. Since $V$ is a faithful $\bar{U}$-module,
it follows from Lemma \ref{lvacuumlikevector} 
and Proposition \ref{pvacuumlikevector2} that the linear map
$f$ from $\bar{U}$ to $V$ such that $f(u)=u_{-1}{\bf 1}$ for $u\in \bar{U}$
is an one-to-one $\bar{U}$-homomorphism.
Clearly, $f$ is onto. Then $V=\bar{U}$ has a 
natural generalized vertex algebra structure.
The other assertions are clear.$\;\;\;\;\Box$

Similar to Lemma \ref{lvacuumlikevector} we have the following result:

\bl{ldderivative}
Let $V$ be a generalized vertex algebra with a graded generating subspace
$A$, $W$ be a $V$-module and $D_{W}$ be a grading-preserving endomorphism 
of $W$ such that
\begin{eqnarray}\label{ederivativeA}
[D_{W},Y(v,z)]={d\over dz}Y(v,z)\;(=Y(D(v),z)) \;\;\;\mbox{ for }v\in A.
\end{eqnarray}
Then (\ref{ederivativeA}) for all $v\in V$.
\el

\pf Recall that $Y(D(v),z)={d\over dz}Y(v,z)$ for all $v\in V$.
Set
\begin{eqnarray}
K=\{ v\in V\;|\; [D_{W},Y(v,z)]=Y(D(v),z)={d\over dz}Y(v,z)\}.
\end{eqnarray}
Then $A\subset K$. 

Let $u\in K\cap V^{g},\; v\in K\cap V^{h},\; g,h\in G$ and 
let $w\in W^{s},\; s\in S$. Then we have
\begin{eqnarray}
& &z_{2}^{-1}\left(\frac{z_{1}-z_{0}}{z_{2}}\right)^{-(g,s)}
\delta\left(\frac{z_{1}-z_{0}}{z_{2}}\right)
[D_{W},Y(Y(u,z_{0})v,z_{2})]w\nonumber\\
&=&z_{0}^{-1}\left(\frac{z_{1}-z_{2}}{z_{0}}\right)^{(g,h)}
\delta\left(\frac{z_{1}-z_{2}}{z_{0}}\right)[D_{W},Y(u,z_{1})]Y(v,z_{2})w
\nonumber\\
& &+z_{0}^{-1}\left(\frac{z_{1}-z_{2}}{z_{0}}\right)^{(g,h)}
\delta\left(\frac{z_{1}-z_{2}}{z_{0}}\right)Y(u,z_{1})[D_{W},Y(v,z_{2})]w
\nonumber\\
& &-c(g,h)z_{0}^{-1}\left(\frac{z_{2}-z_{1}}{z_{0}}\right)^{(g,h)}
\delta\left(\frac{z_{2}-z_{1}}{-z_{0}}\right)[D_{W},Y(v,z_{2})]Y(u,z_{1})w
\nonumber\\
& &-c(g,h)z_{0}^{-1}\left(\frac{z_{2}-z_{1}}{z_{0}}\right)^{(g,h)}
\delta\left(\frac{z_{2}-z_{1}}{-z_{0}}\right)Y(v,z_{2})[D_{W},Y(u,z_{1})]w
\nonumber\\
&=&z_{0}^{-1}\left(\frac{z_{1}-z_{2}}{z_{0}}\right)^{(g,h)}
\delta\left(\frac{z_{1}-z_{2}}{z_{0}}\right)
Y(D(u),z_{1})Y(v,z_{2})w
\nonumber\\
& &+z_{0}^{-1}\left(\frac{z_{1}-z_{2}}{z_{0}}\right)^{(g,h)}
\delta\left(\frac{z_{1}-z_{2}}{z_{0}}\right)
Y(u,z_{1})Y(D(v),z_{2})w\nonumber\\
& &-c(g,h)z_{0}^{-1}\left(\frac{z_{2}-z_{1}}{z_{0}}\right)^{(g,h)}
\delta\left(\frac{z_{2}-z_{1}}{-z_{0}}\right)
Y(D(v),z_{2})Y(u,z_{1})w
\nonumber\\
& &-c(g,h)z_{0}^{-1}\left(\frac{z_{2}-z_{1}}{z_{0}}\right)^{(g,h)}
\delta\left(\frac{z_{2}-z_{1}}{-z_{0}}\right)
Y(v,z_{2})Y(D(u),z_{1})w\nonumber\\
&=&z_{2}^{-1}\left(\frac{z_{1}-z_{0}}{z_{2}}\right)^{-(g,s)}
\delta\left(\frac{z_{1}-z_{0}}{z_{2}}\right)
(Y(Y(D(u),z_{0})v,z_{2})+Y(Y(u,z_{0})D(v),z_{2}))w\nonumber\\
&=&z_{2}^{-1}\left(\frac{z_{1}-z_{0}}{z_{2}}\right)^{-(g,s)}
\delta\left(\frac{z_{1}-z_{0}}{z_{2}}\right)
Y(DY(u,z_{0})v,z_{2})w.
\end{eqnarray}
This gives
\begin{eqnarray}
[D_{W},Y(Y(u,z_{0})v,z_{2})]w=Y(DY(u,z_{0})v,z_{2})w.
\end{eqnarray}
Then $Y(u,z_{0})v\in K\{z_{0}\}$. Thus $K$ is a subalgebra of $V$
containing $A$ and ${\bf 1}$.
Since $A$ generates $V$, we must have $K=V$. This concludes the proof.
$\;\;\;\;\Box$

With Lemma \ref{ldderivative} we immediately have:

\bp{plocalderivative}
Let $A$ be a generalized vertex pre-algebra of parafermions on $W$ and
let $D_{W}$ be a grading-preserving endomorphism of $W$ such that
\begin{eqnarray}
[D_{W},a(z)]={d\over dz}a(z)\;\;\;\mbox{ for }a(z)\in A.
\end{eqnarray}
Denote by $V$ the generalized vertex algebra generated by $A$.
Then on $W$,
\begin{eqnarray}
[D_{W},Y(v,z)]=Y(D(v),z)={d\over dz}Y(v,z)\;\;\;\mbox{ for }v\in V.
\;\;\;\;\Box
\end{eqnarray}
\ep

\newpage

\section{Generalized vertex algebras associated with $Z$-algebras}
In this section we briefly recall from [LW1-3] and [LP1-2]
the fundamental results about $Z$-operators and we then
show that the vacuum space $\Omega_{M(\ell,0)}$  of the 
generalized Verma module $M(\ell,0)$ 
for an affine Lie algebra $\hat{\fg}$ with respect to 
the homogeneous Heisenberg subalgebra has a canonical generalized 
vertex algebra structure. We also show that for any highest weight
$\hat{\fg}$-module $W$ of level $\ell$, $\Omega_{W}$ is an
$\Omega_{M(\ell,0)}$-module.

Let $\fg$ be a finite-dimensional simple Lie algebra, ${\bf h}$ be a
Cartan subalgebra, $\Phi$ be the set of roots and
$\Q$ be the root lattice.
Let $\<\cdot,\cdot\>$ be 
the normalized killing form such that $\<\a,\a\>=2$ for a long root $\a$.
Using $\<\cdot,\cdot\>$ we identify ${\bf h}^{*}$ with ${\bf h}$.

Let $\hat{\fg}$ be the affine Lie algebra:
\begin{eqnarray}
\hat{\fg}=\fg\otimes \C[t,t^{-1}]\oplus\C c,
\end{eqnarray}
where
\begin{eqnarray}	
& &[a\otimes t^{m},b\otimes t^{n}]
=[a,b]\otimes t^{m+n}+m\<a,b\>\delta_{m+n,0}c,\label{eaffinerelation}\\
&&[\hat{\fg},c]=0
\end{eqnarray}
for $a,b\in \fg,\; m,n\in \Z$. Following the tradition, we also
use $a(n)$ for $a\otimes t^{n}$. For $n\in \Z$, we denote
\begin{eqnarray}
\fg (n)=\{ a(n)\;|\; a\in \fg\}.
\end{eqnarray}
For $a\in \fg$, define the generating series
\begin{eqnarray}
a(z)=\sum_{n\in \Z}(a\otimes t^{n})z^{-n-1}\in \hat{\fg}[[z,z^{-1}]].
\end{eqnarray}

Set
\begin{eqnarray}
\hat{\bf h}^{+}={\bf h} \otimes t\C[t],\;\;\;\;
\hat{\bf h}^{-}={\bf h} \otimes t^{-1}\C[t^{-1}],
\end{eqnarray}
subalgebras of $\hat{\bf h}={\bf h}\otimes \C[t,t^{-1}]+\C c$. Then
\begin{eqnarray}
\hat{\bf h}_{\Z}=\hat{\bf h}^{+}+\hat{\bf h}^{-}+\C c
\end{eqnarray}
is a Heisenberg subalgebra of $\hat{\fg}$. 
For $\ell\in \C^{\times}$, 
let $M(\ell)$ be the standard irreducible $\hat{\bf h}_{\Z}$-module 
with $c$ acting as scalar $\ell$. We may also consider
$M(\ell)$ as an $\hat{\bf h}$-module with the action of ${\bf h}$ being zero.

\bd{dcategory2}
{\em For $\ell \in \C$, we define a category ${\cal{C}}_{\ell}$
whose objects are level-$\ell$ $\hat{\fg}$-modules $W$ which are
${\bf h}$-weight modules 
satisfying the condition that for every $w\in W$, 
$\fg (n)w=0$ for $n$ sufficiently large and 
$\dim U(\hat{\bf h}^{+})w<\infty$.}
\ed

It follows from [LW3] and [K] that each $W$ from ${\cal{C}}_{\ell}$ 
with $\ell\ne 0$ is a completely reducible $\hat{\bf h}$-module.
For $W\in {\cal{C}}_{\ell}$, set
\begin{eqnarray}
\Omega_{W}
=\{ w\in W\;|\; h(n)w=0\;\;\;\mbox{ for }h\in {\bf h},\; n\in \Z_{+}\}.
\end{eqnarray}
Since $[h(0), h'(m)]=0$ for $h,h'\in {\bf h},\;m\in \Z$, 
$h(0)$ preserves $\Omega_{W}$.
With $\hat{\fg}$, as an ${\bf h}$-module,
being naturally ${\bf h}\;(={\bf h}^{*})$-graded we have
\begin{eqnarray}
\Omega_{W}=\coprod_{\a\in {\bf h}}\Omega_{W}^{\a},
\end{eqnarray}
where 
$\Omega_{W}^{\a}=\{ w\in \Omega_{W}\;|\; h(0)w=\<\a,h\>w
\;\;\;\mbox{ for }h\in {\bf h}\}$.

For $h\in {\bf h}$, set ([LW1])
\begin{eqnarray}
E^{\pm}(h,z)
=\exp \left(\sum_{\pm n\in \Z_{+}}{h(n)\over n}z^{-n}\right)
\in U(\hat{\fg})[[z^{\mp 1}]].
\end{eqnarray}
For $a\in \fg_{\a},\; \a\in \Phi$, we define
\begin{eqnarray}
Z(a,z)=E^{-}({1\over\ell}\a,z)a(z)E^{+}({1\over\ell}\a,z),
\end{eqnarray}
a formal object. (It is an element of $\overline{U(\hat{\fg})}[[z,z^{-1}]]$,
where $\overline{U(\hat{\fg})}$ is a certain formal completion 
of $U(\hat{\fg})$.) For every $W\in {\cal{C}}_{\ell}$, $Z(v,z)$ is a
well defined element of $(\End W)[[z,z^{-1}]]$.
Then we have ([LW4], [LP2]):

\bp{plwlp1} Let $\ell\in \C^{\times},\; W\in {\cal{C}}_{\ell}$.
For $h\in {\bf h},\; u\in \fg_{\a},\;v\in \fg_{\b},\; \a,\b\in \Phi$, on $W$,
\begin{eqnarray}
& &[h(0),Z(u,z)]=\<\a,h\>Z(u,z),\label{ecompat}\\
& &[h(n),Z(u,z)]=0\;\;\;\mbox{ for }n\ne 0,\label{ecommutting}
\end{eqnarray}
$$(1-z_{2}/z_{1})^{\<\a,\b\>/\ell}Z(u,z_{1})Z(v,z_{2})
-(1-z_{1}/z_{2})^{\<\a,\b\>/\ell}Z(v,z_{2})Z(u,z_{1})=$$
\begin{equation}\label{ezalgebra}
=\left\{
\begin{array}{c}
z_{1}^{-1}\delta(z_{2}/z_{1})Z([u,v],z_{2})\;
\hspace{4cm}\mbox{ if }\a+\b\ne 0,\\
z_{1}^{-1}\delta(z_{2}/z_{1})[u,v]z_{2}^{-1}
+\ell \<u,v\> {\partial \over\partial z_{2}}z_{1}^{-1}\delta(z_{2}/z_{1})
\;\;\;\;\mbox{ if }\a+\b=0.
\end{array}\right.
\end{equation}
\ep

It follows immediately from (\ref{ecommutting}) that
$Z(v,z)$ maps $\Omega_{W}$ to $\Omega_{W}[[z,z^{-1}]]$.

Set
\begin{eqnarray}
Z(v,z)=\sum_{n\in \Z}Z(v,n)z^{-n}.
\end{eqnarray}

\bd{dlp2} {\em 
Following [LP2] we define a category ${\cal{Z}}_{\ell}$ 
(which was denoted by
$P_{\ell}$ in [LP2]) whose objects are ${\bf h}$-weight modules $U$ 
equipped with a family of operators $Z_{U}(a,m)$ (linear in $a$) on $U$ for 
$a\in \fg_{\a},\;\a\in \Phi,\; m\in \Z$ such that 
$Z_{U}(a,z)w\in U((z))$ and such that
(\ref{ecompat}) and 
(\ref{ezalgebra}) hold for $Z_{W}(u,z)$ in place of $Z(u,z)$.}
\ed

Clearly, we have a functor $\Omega$ from ${\cal{C}}_{\ell}$
to ${\cal{Z}}_{\ell}$. Conversely, given $U\in {\cal{Z}}_{\ell}$, we set
\begin{eqnarray}
E(U)=M(\ell)\otimes U\;\;\;\mbox{ as a vector space. }
\end{eqnarray}
View $E(U)$ as a natural $\hat{\bf h}_{\Z}$-module with
$h(n)$ (for $n\ne 0$) acting on the first factor. 
Then define an action of $\hat{\fg}$ by
\begin{eqnarray}
& &c\mapsto \ell,\\
& &h\mapsto 1\otimes h\;\;\;\mbox{ for }h\in {\bf h},\\
& &a(z)\mapsto  E^{-}(-{1\over\ell}\a,z)
E^{+}(-{1\over\ell}\a,z)\otimes Z_{U}(a,z)
\;\;\;\mbox{ for }a\in \fg_{\a},\; \a\in \Phi.
\end{eqnarray}
It was proved in [LP2] that $E(U)$ is a $\hat{\fg}$-module 
in ${\cal{C}}_{\ell}$.
Furthermore, we have ([LW4], [LP2]):

\bp{plwlp2}
Let $\ell\in \C^{\times}$. Then the functors
\begin{eqnarray}
\Omega: W\mapsto \Omega_{W}\;\;\mbox{and }\;\; E: U\mapsto E(U)
\end{eqnarray}
are exact and they define equivalences
 between the categories ${\cal{C}}_{\ell}$ and
${\cal{Z}}_{\ell}$.
In particular,
$W$ is irreducible in ${\cal{C}}_{\ell}$ if and only if 
$\Omega_{W}$ is irreducible in ${\cal{Z}}_{\ell}$.
\ep

For $\lambda\in {\bf h}$, let $M(\ell,\lambda)$ be 
the Verma $\hat{\fg}$-module. In view of 
the universal property for $M(\ell,\lambda)$,
with Proposition \ref{plwlp2} we immediately have:

\bc{cuniversal}
Let $\ell\in \C,\; \lambda\in {\bf h}$ and let $v$ be a
(nonzero) highest weight vector of $M(\ell,\lambda)$. Let
$U\in {\cal{Z}}_{\ell}$ and let $e\in U$ satisfying
the following conditions:
\begin{eqnarray}
& &he=\<h,\lambda\>e\;\;\;\mbox{ for }h\in {\bf h},\\
& &Z_{U}(u,z)e\in z^{-1}U[[z]]\;\;\;\mbox{ for }u\in \fg_{\a},\; \a\in \Phi,\\
& &Z_{U}(v,z)e\in U[[z]]\;\;\;\mbox{ for }v\in \fg_{\b},\; \b\in \Phi_{+}.
\end{eqnarray}
Then there exists a unique morphism in ${\cal{Z}}_{\ell}$ from 
$\Omega_{M(\ell,\lambda)}$ to $U$ sending $v$ to $e$. $\;\;\;\;\Box$
\ec



\bd{dpsioperator}
{\em For $a\in \fg_{\a},\; \a\in \Phi$, we define
\begin{eqnarray}
\psi (a, z)=Z(a,z)z^{-{1\over\ell}\a(0)}
=E^{-}({1\over\ell}\a,z)a(z)E^{+}({1\over\ell}\a,z)z^{-{1\over\ell}\a (0)}.
\end{eqnarray}}
\ed
Then
\begin{eqnarray}
[h(n),\psi (a,z)]=0\;\;\;\mbox{ for }h\in {\bf h},\; n\ne 0,
\end{eqnarray}
hence $\psi (a,z)$ maps $\Omega_{W}$ to $\Omega_{W}\{z\}$ 
for $W\in {\cal{C}}_{\ell}$. 
Note that (\ref{ecompat}) amounts to
\begin{eqnarray}\label{ecommhz}
z^{h(0)}Z(u,z_{1})=Z(u,z_{1})z^{\<\a,h\>+h(0)}
\end{eqnarray}
for $h\in {\bf h},\; u\in \fg_{\a},\; \a\in \Phi$.
We have the following reformulation of
Proposition \ref{plwlp1} in terms of $\psi$-operators:

\bp{plp2consequence}
Let $u\in \fg_{\a},\; v\in \fg_{\b},\;\a,\b\in \Phi,\; \ell \in \C^{\times}$ 
and $W\in {\cal{Z}}_{\ell}$. On $W$, 
\begin{eqnarray}
& &[h(0),\psi(u,z)]=\<\a,h\>\psi (u,z),\label{ehpsi}\\
& &[h(n),\psi(u,z)]=0\;\;\;\mbox{ for }h\in {\bf h},\;n\ne 0,
\end{eqnarray}
$$(z_{1}-z_{2})^{\<\a,\b\>/\ell}\psi(u,z_{1})\psi(v,z_{2})
-(z_{2}-z_{1})^{\<\a,\b\>/\ell}\psi(v,z_{2})\psi(u,z_{1})=$$
\begin{eqnarray}\label{epsigcomm}
=\left\{\!\!
\begin{array}{c}
z_{1}^{-1}\delta(z_{2}/z_{1})\psi([u,v],z_{2})
(z_{2}/z_{1})^{{1\over\ell}\a(0)}
\hspace{1cm}\mbox{ if }\a+\b\ne 0,\\
\ell \<u,v\> 
{\partial\over\partial z_{2}}\left(z_{1}^{-1}\delta(z_{2}/z_{1})
(z_{2}/z_{1})^{{1\over \ell}\a(0)}\right)\;\;\;\;\;\;\mbox{ if }\a+\b=0.
\end{array}\right.
\end{eqnarray}
\ep

\pf The first two identities are obvious. Using (\ref{ecommhz}) we obtain
$$(z_{1}-z_{2})^{\<\a,\b\>/\ell}\psi(u,z_{1})\psi(v,z_{2})
-(z_{2}-z_{1})^{\<\a,\b\>/\ell}\psi(v,z_{2})\psi(u,z_{1})=$$
\begin{eqnarray}
=\left\{\!\!
\begin{array}{c}
z_{1}^{-1}\delta(z_{2}/z_{1})\psi([u,v],z_{2})
(z_{2}/z_{1})^{{1\over\ell}\a(0)}\;
\hspace{4.5cm}\mbox{ if }\a+\b\ne 0,\\
\left(z_{1}^{-1}\delta(z_{2}/z_{1})[u,v]z_{2}^{-1}
+\ell \<u,v\> {\partial\over\partial z_{2}}z_{1}^{-1}\delta(z_{2}/z_{1})\right)
(z_{2}/z_{1})^{{1\over \ell}\a(0)}\;\;\;\;\mbox{ if }\a+\b=0.
\end{array}\right.
\end{eqnarray}
It remains to consider the case $\a+\b=0$.
Using the fact that $[u,v]=\<u,v\>\a$ and $\delta(z)z^{m}=\delta(z)$ 
for $m\in \Z$ we obtain
\begin{eqnarray}
& &\left(z_{1}^{-1}\delta(z_{2}/z_{1})[u,v]z_{2}^{-1}
+\ell \<u,v\> {\partial\over\partial z_{2}}z_{1}^{-1}\delta(z_{2}/z_{1})\right)
(z_{2}/z_{1})^{{1\over \ell}\a(0)}\nonumber\\
&=&\ell \<u,v\> {\partial\over\partial z_{2}}
\left(z_{1}^{-1}\delta(z_{2}/z_{1})(z_{2}/z_{1})^{{1\over\ell}\a(0)}\right).
\end{eqnarray}
This completes the proof.$\;\;\;\;\Box$

It is a simple fact (see for example [Li2]) that
\begin{eqnarray}\label{ezerodelta}
(z_{1}-z_{2})^{m}\left({\partial\over\partial z_{2}}\right)^{n}
z_{2}^{-1}\delta(z_{1}/z_{2})=0
\end{eqnarray}
for $m,n\in \Z$ with $m>n\ge 0$. Then using 
Proposition \ref{plp2consequence} we get
\begin{eqnarray}
(z_{1}-z_{2})^{2}\left((z_{1}-z_{2})^{\<\a,\b\>/\ell}\psi(u,z_{1})\psi(v,z_{2})
-(z_{2}-z_{1})^{\<\a,\b\>/\ell}\psi(v,z_{2})\psi(u,z_{1})\right)=0,
\end{eqnarray}
where $u,v$ are as in Proposition \ref{plp2consequence}.
Then for any $U\in {\cal{Z}}_{\ell}$, e.g., 
$U=\Omega_{W}$ for some $W\in {\cal{C}}_{\ell}$,
$\psi(u,z)$ for $u\in \fg_{\a},\;\a\in \Phi$ linearly span
a generalized vertex pre-algebra of parafermion operators on $U$, which
by Theorem 3.16 generates a canonical generalized vertex algebra
inside $(\End U)\{z\}$ with $G={\bf h},\; c(\cdot,\cdot)=1,\; 
(\cdot,\cdot)={1\over\ell}\<\cdot,\cdot\>$.

\bl{ladjointmodule}
Let $U\in {\cal{Z}}_{\ell}$ and let $V$ be the
generalized vertex algebra generated by $\psi(u,z)$ for 
$u\in \fg_{\a},\;\a\in \Phi$ inside $(\End U)\{z\}$.
Then $V$ is a natural object of ${\cal{Z}}_{\ell}$ where
\begin{eqnarray}
& &h\cdot \phi(z)=[h,\phi(z)],\\
& &Z_{V}(u,z_{0})=Y_{V}(\psi(u,z),z_{0})z_{0}^{{1\over\ell}\a(0)}
\end{eqnarray}
for $h\in {\bf h},\;\phi(z)\in V,\; u\in \fg_{\a},\; \a\in \Phi$.
\el

\pf Since $U$ is an ${\bf h}$-weight module, 
$\End U$ is a natural ${\bf h}$-module where
\begin{eqnarray}
h\cdot f=[h,f]\;(=hf-fh)\;\;\;\mbox{ for }h\in {\bf h},\;f\in \End U.
\end{eqnarray}
Then $(\End U)\{z\}$ is a natural ${\bf h}$-module.
Since the generators
$\psi(a,z)$ for $a\in \fg_{\a}$ of $V$ are ${\bf h}$-eigenvectors
(recall (\ref{ehpsi})),
using the proof of Lemma \ref{ldderivative} we can easily
show that $V$ is an ${\bf h}$-weight module and (\ref{ehpsi}) holds on $V$. 
Note that $U$ is a faithful $V$-module.
Then it follows immediately from Lemma \ref{lequivalence} that
(\ref{epsigcomm}) holds on $V$. This shows that $V$ is a natural 
object of ${\cal{Z}}_{\ell}$. $\;\;\;\;\Box$

Let $\ell\in \C^{\times}$. Consider 
the generalized Verma $\hat{\fg}$-module
\begin{eqnarray}
M(\ell,0)=U(\hat{\fg})\otimes_{U(\fg\otimes \C[t]+\C c)}\C_{\ell},
\end{eqnarray}
where $\C_{\ell}=\C$ as a vector space and $\fg\otimes \C[t]$ 
acts as zero on $\C_{\ell}$ and $c$ acts as $\ell$. Denote by ${\bf 1}$
the highest weight vector $1\otimes 1$ of $M(\ell,0)$.
Let $L(\ell,0)$ be the (unique) irreducible quotient module 
with ${\bf 1}$ as a fixed highest weight vector.
Identify $\fg$ as a subspace of $M(\ell,0)$ 
and $L(\ell,0)$ through $a\mapsto a(-1){\bf 1}$.
Then we have
\begin{eqnarray}
\fg_{\a}\subset \Omega_{M(\ell,0)}^{\a}\subset \Omega_{M(\ell,0)}
\;\;\;\mbox{ for }\a\in \Phi.
\end{eqnarray}

\bt{tcosetalgebraaffine}
Let $\ell\in \C^{\times}$ and $V=M(\ell,0)$ or $L(\ell,0)$. 
Then there exists 
a unique generalized vertex algebra structure $Y_{\Omega}$ 
on $\Omega_{V}$ with $G=\Q, \; c(\cdot,\cdot)=1$
and $(\cdot,\cdot)=\<\cdot,\cdot\>/\ell$ such that
$Y_{\Omega}({\bf 1},z)=1$ and
$Y_{\Omega}(a,z)=\psi(a,z)$ for $u\in \fg_{\a},\; \a\in \Phi$.
Furthermore, $\Omega_{V}$ 
is generated by $\fg_{\a}$ ($\a\in \Phi$).
\et

\pf Clearly, $V$ is $\Q$-graded. Then we take $G=\Q$,
$(\cdot,\cdot)=\<\cdot,\cdot\>/\ell$ and $c(\cdot,\cdot)=1$.
Let $A$ be the linear span of $\psi(a,z)$ for 
$a\in \fg_{\a},\; \a\in \Phi$.
It follows from Proposition \ref{plp2consequence} and 
(\ref{ezerodelta}) that
$A$ is a generalized vertex pre-algebra. 
Since $\Omega_{V}$
is generated from ${\bf 1}$ by all the components of 
$Z(a,z)=\psi(a,z)z^{\a(0)/\ell}$,
there exists a unique generalized vertex algebra structure 
$Y_{\Omega}$ on $\Omega_{V}$ with the required 
conditions.$\;\;\;\;\Box$

\bp{pvoagvoa}
Let $\ell\in \C^{\times}$ and let 
$U\in {\cal{Z}}_{\ell}$. Then $U$ is a natural
$\Omega_{M(\ell,0)}$-module. In particular,
$\Omega_{W}$ is an
$\Omega_{M(\ell,0)}$-module for $W\in {\cal{C}}_{\ell}$.
\ep

\pf Set $M=\Omega_{M(\ell,0)}\oplus U$, an object of ${\cal{Z}}_{\ell}$.
Let $V$ be the generalized vertex algebra generated by 
$\psi(v,z)$ for $v\in \fg_{\a},\; \a\in \Phi$ inside
$(\End \Omega_{M(\ell,0)})\{z\}\subset (\End M)\{z\}.$ 
Then $M$ is a $V$-module
with $\Omega_{M(\ell,0)}$ and $U$ as submodules.
{}From Proposition 3.19, there is a $V$-homomorphism $f$ from
$V$ onto $\Omega_{M(\ell,0)}$, which maps $I(z)$ to ${\bf 1}$. 
In view of Lemma \ref{ladjointmodule}, $V$ is a natural
${\bf h}$-module in ${\cal{Z}}_{\ell}$.
It follows from Corollary 4.5 that
$f$ is a linear isomorphism, hence $V=\Omega_{M(\ell,0)}$.
Thus $M$ is an $\Omega_{M(\ell,0)}$-module. Therefore, 
$U$ is an $\Omega_{M(\ell,0)}$-module.$\;\;\;\;\Box$

We define an $\Omega_{M(\ell,0)}$-${\bf h}$-module to be an 
$\Omega_{M(\ell,0)}$-module and an ${\bf h}$-weight module such that
(\ref{ehpsi}) holds. In view of Proposition \ref{pvoagvoa} and Lemma 
\ref{ladjointmodule} we immediately have:

\bc{ccategoryequivalence}
The category ${\cal{Z}}_{\ell}$ is canonically isomorphic to
the category of $\Omega_{M(\ell,0)}$-${\bf h}$-modules.
$\;\;\;\;\Box$
\ec

\br{rgvoadl2}
{\em Note that for a positive integer $\ell$, the generalized 
vertex operator algebra $\Omega_{L(\ell,0)}^{B}$ 
constructed in [DL2] is a quotient algebra of $\Omega_{M(\ell,0)}$.
This is recently studied in [Li4] from a different point of view. }
\er

\end{document}